\documentclass[11pt]{article}
\usepackage[centertags,intlimits]{amsmath}
\usepackage{amsfonts}
\usepackage{amsthm}
\usepackage{graphics}
\allowdisplaybreaks[2]
\numberwithin{equation}{section}
\newtheorem{theorem}{Theorem}[section]
\newtheorem{lemma}{Lemma}[section]

\theoremstyle{remark}
\newtheorem{remark}{Remark}[section]

\providecommand{\abs}[1]{\lvert #1\rvert}

\newcommand{\nc}{\newcommand}
\nc{\vb}{\mathbf{v}}
\nc{\bx}{\mathbf{x}}
\nc{\by}{\mathbf{y}}
\nc{\bz}{\mathbf{z}}
\nc{\bu}{\mathbf{u}}
\nc{\bv}{\mathbf{v}}
\nc{\ba}{\mathbf{a}}
\nc{\bs}{\mathbf{s}}
\nc{\bq}{\mathbf{q}}
\nc{\bd}{\mathbf{d}}
\nc{\bb}{\mathbf{b}}
\nc{\bc}{\mathbf{c}}
\nc{\bi}{\mathbf{i}}
\nc{\bfr}{\mathbf{r}}
\nc{\bP}{\mathbf{P}}
\nc{\bQ}{\mathbf{Q}}
\nc{\R}{\mathbb R}
\nc{\N}{\mathbb N}
\nc{\bbC}{\mathbb C}
\nc{\D}{\mathbb D}
\nc{\Z}{\mathbb Z}
\nc{\F}{\mathbf F}
\nc{\bbS}{\mathbb S}
\nc{\bE}{\mathbf E}
\nc{\B}{\cal B}
\nc{\br}{\bigr}
\nc{\bl}{\bigl}
\nc{\Bl}{\Bigl}
\nc{\Br}{\Bigr}
\nc{\ind}{\mathbf{1}}

\title{On the Limit Law of a Random Walk\\ Conditioned  to Reach a High Level}
\author{Sergey G. Foss\footnote{Email address: S.Foss@ma.hw.ac.uk} \\
 Heriot-Watt University, Edinburgh, UK and\\
Institute of Mathematics, Novosibirsk, Russia\\ and\\
Anatolii A. Puhalskii\footnote{Email address: anatolii.puhalskii@ucdenver.edu}\\
University of Colorado Denver, Denver, U.S.A. and\\
 Institute for Problems in Information
Transmission, Moscow, Russia}
\begin{document}
\maketitle
\sloppy
\begin{abstract}
We consider a random walk with a negative drift and with a jump distribution
which under Cram\'er's
change of measure belongs  to the domain of attraction of a
spectrally positive stable law. If conditioned to reach a
high level and suitably scaled, this random walk converges in law to a nondecreasing Markov
process which can be interpreted as a spectrally-positive L\'evy 
process conditioned not to
overshoot level one.

\end{abstract}

\section{Introduction}
\label{sec:overview}

Let $\xi_1,\xi_2,\ldots$ be i.i.d. random variables on a probability space $(\Omega,\mathcal{F},\bP)$
 with $\bE\xi_1<0$ and $\mathbf{P}(\xi_1>0)>0$. Then the random walk
 $S_n=\sum_{i=1}^n\xi_i$ tends to $-\infty$ with probability 1 as
 $n\to\infty$ and the
 event that it exceeds a high  level has a small albeit
 positive probability.
The asymptotics of this probability
 have been studied extensively.
Let ${\mathbf E} \exp (\lambda \xi_1)<\infty$ for some
$\lambda>0$ and let
\begin{equation}
  \label{eq:43}
\gamma = \sup \{ \lambda \ : \ {\mathbf E} \exp (\lambda \xi_1) \le 1 \}.
\end{equation}
Clearly,   $\gamma>0$ and
 $\mathbf{E}\exp(\gamma\xi_1)\le1$\,.

In the   ``classical'' case where
 $\mathbf{E}\exp(\gamma\xi_1)=1$ and
$\beta=\mathbf{E}\,\xi_1\exp(\gamma\xi_1)<\infty$,
the celebrated Cram\'er-Lundberg
theorem asserts that, for a certain constant $C_1$,
\begin{equation}\label{classic}
{\mathbf P} (\sup_{n}S_n>r) \sim C_1 e^{-\gamma r} \quad\text{ as } r\to\infty,
\end{equation}
where $\sim$ stands for asymptotic equivalence.
The limit is taken along
all $r$ if  $\xi_1$ has a non-lattice distribution
and along multiples of the lattice span if  $\xi_1$ has a lattice
distribution, see, for example, Asmussen \cite[XIII.5]{Asm03},
Borovkov \cite[\S22]{Bor72}, or Feller \cite[XII]{Fel71}.

If $\mathbf{E}\exp(\gamma\xi_1)<1$ so that
$\mathbf{E}\exp(\lambda\xi_1)=\infty$ for all $\lambda>\gamma$, then,
under certain regularity
assumptions on the distribution of  $\xi_1$ (more specifically,
provided it belongs to  class ${\cal S}_{\gamma}$, see
Teugels \cite{Teu75}),
\begin{equation*}
{\mathbf P}(\sup_{n}S_n>r) \sim C_2 {\mathbf P} (\xi_1>r)\;\; \text{ as }r\to\infty,
\end{equation*}
where
$C_2 =
 {\mathbf E} \exp (\gamma \sup_{n}S_n)/\bl(1-\mathbf{E}\xi_1\exp(\gamma\xi_1)
 \br)$, see, for example, Bertoin and Doney \cite{BerDon} and references
 therein. For earlier results, see Borovkov \cite[\S22]{Bor72};
 recent developments can be found in Borovkov and Borovkov
\cite{BorBor04} and Zachary and Foss \cite{ZF}.

The     borderline case where
$\mathbf{E}\exp(\gamma\xi_1)=1$ and
$\beta=\infty$ was first addressed in Borovkov's monograph, Borovkov
\cite[\S22]{Bor72}. More complete results have been obtained by
Korshunov \cite{Kor05} who showed   that if the  distribution of
 $\xi_1$ is nonlattice and
 the distribution function $F$, where     $F(dx)= \exp(\gamma x)\mathbf{P}(\xi_1\in dx)$,
has a
regularly varying right-hand tail with index $-\alpha\in (-1,-1/2)$, then
\begin{equation} \label{Dima}
{\mathbf P}(\sup_{n}S_n>r) \sim C_3\,
\frac{e^{-\gamma r}}{\gamma m(r)}\;\;\text{ as }r\to\infty,
\end{equation}
for a certain universal constant $C_3$ and for $m(r) = \int_0^r \,
(1-F(u )) \, du$\,. A similar asymptotic relationship was
shown to be valid for lattice distributions too, with a different
constant $C_3^{'}$. 

One can also clarify in what way
 the  event of attaining a high level $r$ is most likely to occur.
In the classical case, the trajectory $(S_{\lfloor
 rt\rfloor}/r,\,t\in\R_+)$, where $\lfloor \cdot\rfloor$ stands for the
 integer part,
 stays with a probability close to one  in a neighbourhood
 of the   straight line with slope $\beta$;   see  (\ref{eq:70})
below for the  exact formulation.
Thus, the event of reaching
 level $r$ is  realised typically via multiple (of order $r$) jumps of order 1.
If ${\mathbf E} \exp (\gamma \xi_1) <1$, then the high level
is most likely to be reached at the very beginning of the random
walk, which occurs due to a single big jump of order $r$, see Borovkov
and Borovkov \cite{BorBor04} for the case of regular exponential
distribution tails and Zachary and Foss \cite{ZF} for the more general
${\cal S}_{\gamma}$ distributions.
In fact, as it  follows from the results of
Zachary and Foss \cite{ZF},
the conditional distribution of the  jump time
converges weakly (no  scaling is involved) to the geometric
distribution with parameter $p=1-{\mathbf E} \exp (\gamma \xi_1)$\,.

The purpose of this paper is to  study
 the most likely way for the random walk to attain
a high level $r$ in the setting considered by Korshunov  \cite{Kor05}.
Not unexpectedly,
the results and intuitive explanations in the borderline situation
  are essentially more intriguing and
complicated.
On the one hand, as $\beta\uparrow \infty$ in the classical case,
 the trajectories $S_{\lfloor rt\rfloor}/r$ for large values of $r$
increasingly look like a jump at time $0$. On the other hand,
 as ${\mathbf E} \exp (\gamma \xi_1)\uparrow 1$, the time for the big
 jump that reaches (or almost reaches) $r$ to occur tends to infinity.
Therefore, the typical jump sizes in the case we concern ourselves with
 here should, on the one hand, grow to infinity as $r\to\infty$, but, on the
    other hand, be  of a smaller order of magnitude than $r$.
The correct order of magnitude is captured by considering the random
 walk $S_n$ with $n$ tending to infinity at a slower rate than
 $r$\,. Typically, one takes
 $n=\lfloor r^\alpha t\rfloor$\,, which corresponds to the 
  jumps  of order
 $r^{1-\alpha}$\,. 

 To  provide better insight into the kind of results we obtain,
 let us recall the  argument
 underlying the asymptotics  in the classical case. Its main points can
 be found in  Asmussen \cite{Asm82}, who refers to Iglehart
 \cite{Igl72} and von  Bahr \cite{Bah74}.
Let $\mathcal{F}_n$ denote the $\sigma$-algebra on $\Omega$ generated
by the $\xi_i,\,i=1,2,\ldots,n$. We shall assume that
the $\sigma$-algebra $\mathcal{F}$ is generated by the
$\sigma$-algebras $\mathcal{F}_n,\,n=1,2,\ldots$\,.
Let   $\bP^\ast$ be the measure  on $(\Omega,\mathcal{F})$ 
defined by
\begin{equation}
  \label{eq:48}
\bP^\ast(\Gamma)=\bE\exp(\gamma
S_n)\ind_\Gamma\;\; \text{ for } \Gamma\in\mathcal{F}_n,
\end{equation}
where   $\ind_\Gamma$ denotes the indicator function of event $\Gamma$\,.
 It is a probability
measure by the assumption that $\mathbf{E}\exp(\gamma\xi_1)=1$\,.  The probability
measures $\bP$ and $\bP^\ast$ are locally equivalent and
 $d\bP/d\bP^\ast\big|_{\mathcal{F}_n}=\exp(-\gamma S_n)$.
We
also note that under $\bP^\ast$ the $\xi_k$ are i.i.d. with
mean $\beta$\,.

For $r>0$, let $\tau^{(r)}$ be the first time the random walk $S_n$
attains level $r$, i.e.,
\begin{equation}
  \label{eq:65}
  \tau^{(r)}=\min\{n:\,S_n\ge r\}\,.
\end{equation}
Since $\{\tau^{(r)}=n\}\in\mathcal{F}_n$,
\begin{equation*}
  \bP(\tau^{(r)}<\infty)=\sum_{n=1}^\infty \bP(\tau^{(r)}=n)=
\sum_{n=1}^\infty \bE^\ast e^{-\gamma S_n}\ind_{\{\tau^{(r)}=n\}}=
\bE^\ast e^{-\gamma S_{\tau^{(r)}}}\ind_{\{\tau^{(r)}<\infty\}}\,,
\end{equation*}
where $\mathbf{E}^\ast$ denotes expectation with respect to
$\mathbf{P}^\ast$\,.
On noting that
 $\bP^\ast(\tau^{(r)}<\infty)=1$ as $  \bE^\ast \xi_1>0$,
we conclude that
\begin{equation}
  \label{eq:12}
  \bP(\tau^{(r)}<\infty)=\bE^\ast \exp(-\gamma S_{\tau^{(r)}}).
\end{equation}

More generally, if $\Gamma\in\mathcal{F}_{\tau^{(r)}}$,
$\mathcal{F}_{\tau^{(r)}}$ being the $\sigma$-algebra associated with
the stopping time $\tau^{(r)}$,
then by the fact that
$\{\tau^{(r)}=n\}\cap\Gamma\in\mathcal{F}_n$
\begin{equation*}
  \bP\bl(\Gamma\cap\{\tau^{(r)}<\infty\}\br)=
\sum_{n=1}^\infty \bP\bl(\Gamma\cap\{\tau^{(r)}=n\}\br)=
\sum_{n=1}^\infty \bE^\ast e^{-\gamma S_n}\ind_{\Gamma\cap\{\tau^{(r)}=n\}}=
\bE^\ast e^{-\gamma S_{\tau^{(r)}}}\ind_{\Gamma},
\end{equation*}
so
\begin{equation}
  \label{eq:66}
  \bP(\Gamma|\tau^{(r)}<\infty)=
\frac{\bE^\ast e^{-\gamma S_{\tau^{(r)}}}\ind_\Gamma}{\bE^\ast e^{-\gamma S_{\tau^{(r)}}}}=
\frac{\bE^\ast e^{-\gamma \chi^{(r)}}\ind_\Gamma}{\bE^\ast e^{-\gamma \chi^{(r)}}}\,,
\end{equation}
where
\begin{equation}
  \label{eq:69}
    \chi^{(r)}=S_{\tau^{(r)}}-r
\end{equation}
is the overshoot of the random walk $S_n$ over level $r$\,.

Suppose now that $\Gamma$ is the event  $\{\sup_{n\le \tau^{(r)}}
\abs{S_n-\beta n}<\varepsilon r\}$\,, where $\varepsilon>0$ is given.
Since $\beta <\infty$, 
the $\mathbf{P}^\ast$-probability of this event tends to $1$ as $r\to\infty$,
by the strong law of large numbers.
Also, the condition $\beta <\infty$ implies, 
 that the $\chi^{(r)}$ under $\mathbf{P}^\ast$  tend in
distribution, as $r\to\infty$, to a proper  random variable, say
$\chi^{(\infty)}$\,,
provided the distribution
of $\xi_1$  is nonlattice --
see, e.g., Asmussen \cite[VIII.2]{Asm03},
Gut \cite[III.10]{Gut88}, or Feller \cite[XI.4]{Fel71}\,.
Therefore, the quantities $\bE^\ast \exp(-\gamma \chi^{(r)})$ converge to a
positive limit  as $r\to\infty$ and, by \eqref{eq:66},
\begin{equation}
  \label{eq:70}
    \lim_{r\to\infty}\bP(\sup_{n\le \tau^{(r)}}
\abs{S_n-\beta n}<\varepsilon r|\tau^{(r)}<\infty)=1\,.
\end{equation}
(In particular,  \eqref{classic}  with $C_1= {\mathbf E}
e^{-\gamma \chi^{(\infty)}}$ follows by
\eqref{eq:12} and \eqref{eq:69}.) 
This argument breaks down in two places if $\beta =\infty$: we can
no longer rely on the law of large numbers for the random walk and
the $\chi^{(r)}$ might  converge to infinity as $r\to\infty$\,. 

In order to tackle these difficilties,
we will assume, as in Korshunov \cite{Kor05}, that  the distribution
function $F$
  has a regularly varying  right-hand tail with index $-\alpha$, where
  $\alpha\in(1/2,1)$.
We scale ``time''   by
$(1- F (r))^{-1}$
  so that the
 processes $(S_{\lfloor (1- F (r))^{-1} t\rfloor}/r,\,t\in\R_+)$
under $\mathbf{P}^\ast$ converge in
    distribution as $r\to\infty$ to a stable subordinator
    $X=(X(t),\,t\in\R_+)$ with L\'evy measure
    $ \alpha x^{-\alpha-1}\,dx$\,, see Resnick \cite{Res86} or Lemma~\ref{le:conv_levy}.
Under the stated
     assumptions, the random variables $\chi^{(r)}/r$
    converge in distribution to a proper random variable $\chi$, which
    assumes values in $(0,1)$ and has density
    $p_\alpha(x)=(\sin\pi\alpha/\pi)x^{-\alpha}(1+x)^{-1}$. (See
Dynkin \cite[Theorem 2]{Dyn55}, or     Feller \cite[XIV.3]{Fel71}, for the case
    of renewal processes, Sinai \cite{Sin57} for the case of the sums of random
    variables with a stable distribution, the  case in question follows by an application
    of Lemma 2 in Korshunov \cite{Kor05}. ) It is then plausible that in \eqref{eq:66}
    one should be able to replace $\chi^{(r)}$ with $r\chi$ so that
    $\exp(-\gamma\chi^{(r)})$ can be replaced with
$ \exp(-r\gamma \chi)$\,. For large values of $r$, the bulk of the
contribution to $\mathbf{E}^\ast \exp(-r\gamma \chi)$
 comes from the small values of $\chi$, so the right-hand side of
 \eqref{eq:66} should be asymptotically equivalent to $\mathbf{P}^\ast(\Gamma|\chi=0)$\,.
If $\Gamma$ is an event associated with the process
$(S_{\lfloor (1- F (r))^{-1} t\rfloor}/r,\,t\in\R_+)$, then it
should translate in the limit into a similar event associated with
$X$. Besides, since $\chi^{(r)}$ is the overshoot over level
$r$ by $S_n$, we have that $\chi^{(r)}/r$ is the overshoot over level
$1$ of the process $(S_{\lfloor (1- F (r))^{-1}
  t\rfloor}/r,\,t\in\R_+)$\,. That this  process converges to $X$
suggests the conjecture
that $\chi$ should be the overshoot of $X$ above level
$1$.

It thus appears as though 
the conditional distribution of the process
$(S_{\lfloor (1- F (r))^{-1}   t\rfloor\wedge \tau^{(r)}}/r,\,t\in\R_+)$
given that $\tau^{(r)}<\infty$ should  converge  to the conditional
distribution of the process $(X(t\wedge \tau)\,,t\in\R_+)$ given the
event $X(\tau)=1$, where
\begin{equation}
  \label{eq:4}
\tau=\inf\{t:\,X(t)\ge 1\}\,.
\end{equation}
The main result of the paper (Theorem 2.2) confirms this conjecture.
As a consequence, we have that
 if the distribution function $ F $ decays as  $x^{-\alpha}$,
then,
assuming $X$ is defined on
a probability space $(\Omega',\mathcal{F}',\mathbf{P}')$, for $B>0$
and $\varepsilon>0$,
\begin{equation*}
  \lim_{r\to\infty}\mathbf{P}(\sup_{n\le
  \tau^{(r)}}\abs{S_n-Bn^{1/\alpha}}<\varepsilon r|\tau^{(r)}<\infty)=
\mathbf{P}'(\sup_{t\le \tau}\abs{X(t)-Bt^{1/\alpha}}<\varepsilon|X(\tau)=1)\,,
\end{equation*}
which can be regarded as a counterpart of \eqref{eq:70}.

We denote the process $ X$ ``conditioned not to overshoot level 1
one'' by
$\widetilde{X}$. There are various ways to define this process
rigourously.
It  can be obtained by an application of the Doob
$h$-transform associated with the function $h(x)=(1-x)^{\alpha-1}$ to the
process $X$. It is also  the limit of conditional distributions of the processes
$(X(t\wedge \tau)\,,t\in\R_+)$ given
$X(\tau)\le 1+\varepsilon$ as $\varepsilon\to0$\,.
We define it by the property that
 the process $(1-\widetilde{X}(t),\,t\in\R_+)$ is a positive 
selfsimilar Markov
process, as  in Lamperti \cite{Lamp72}.

The main statement on the convergence of the  conditional laws of
$(S_{\lfloor (1- F (r))^{-1}   t\rfloor\wedge \tau^{(r)}}/r,\,t\in\R_+)$
given $\tau^{(r)}<\infty$,  to the law of
$\widetilde X$ as
 $r\to\infty$ is formulated in 
Section~\ref{sec:conv-cond-rand} and then proved 
in Section~\ref{sec:proofs}. 
For the latter, we apply the general results on weak convergence of
semimartingales as in Jacod and Shiryaev \cite{JacShi03} and
Liptser and Shiryayev \cite{lipshir}. A key  is to  compute the
predictable jump measures of the processes in question under the
``conditional'' measures. This is done by applying the results on
the transformations of the predictable triplets under absolutely
continuous changes of the measure,  as in Jacod and Shiryaev \cite{JacShi03} and
Liptser and Shiryayev \cite{lipshir}.
Since  the approaches of Korshunov \cite{Kor05} play an important part,
we also have to   treat the nonlattice and lattice cases separately.
We also provide a representation for the limit process
$\widetilde{X}$ which complements the Lamperti representation.
Auxiliary and supplemental results are collected in the Appendix. 
In Subsection A.1 we present a more complete version of the proof of two
theorems from Korshunov \cite{Kor05} which are important for the
developments in this paper,  Subsection A.2 contains 
some  properties of slowly and regularly varying functions, and
 Subsection A.3
contains a  proof of the
convergence in distribution of the processes
$(S_{\lfloor (1- F (r))^{-1}
   t\rfloor}/r,\,t\in\R_+)$ to $X$ under measure $\mathbf{P}^\ast$
 based on the semimartingale weak convergence theory\,.

We conclude the Introduction with a list of notation and conventions
adopted in the paper. We let $\N$ denote the set of positive integers,
 $\R$ the set of real numbers,
$\mathcal{B}(\R)$ the Borel $\sigma$-algebra on $\R$, and $\R_+$ 
the subset of $\R$ of nonnegative reals. For $x\in\R$ and $y\in\R$, $x\wedge
y=\min(x,y)$, $x\vee y=\max(x,y)$, and $x^-=-x\wedge0$\,. Recall also that
$\lfloor x\rfloor$ denotes the integer part of $x$ and $\ind_\Gamma$
the indicator-function of event $\Gamma$\,.
Two positive functions $f(x)$ and $g(x)$ of a real-valued
argument are said to be asymptotically equivalent as $x\to\infty$, which is written
as $f(x)\sim g(x)$ if $\lim_{x\to\infty}f(x)/g(x)=1$\,. We write
$f(x)=O(g(x))$ if $f(x)\le Cg(x)$ for some $C>0$ and for all $x$ great enough.
Integrals of the form $\int_a^b$ are understood as $\int_{(a,b]}$
unless otherwise indicated. 

We denote  by $\D$
the  space of
real-valued right-continuous   functions on
$\R_+$ with left-hand limits. Its elements
 are denoted by lower-case bold-face Roman characters,
e.g., $\bx=(\bx(t),\,t\in\R_+)$; $\bx(t-)$ is the left-hand
limit of $\bx$ at $t$\, and $\Delta\bx(t)=\bx(t)-\bx(t-)$ is the size of
the jump at $t$\,. The space $\D$ is assumed to be 
 endowed with the Skorohod
 $J_1$-topology, equipped with the Borel $\sigma$-algebra $\mathcal{B}(\D)$,
 and metrised by a complete separable metric,  see
 Ethier and Kurtz \cite{EthKur86},
 Jacod and Shiryaev
\cite{JacShi03}, and Liptser and Shiryaev \cite{lipshir} for the definition and
properties;
$\D_\uparrow$ denotes the subset of
$\D$ of increasing functions with $\bx(0)=0$
which is  equipped  with the subspace topology.
All stochastic processes encountered in this paper have trajectories
  in $\D$ and are considered as random elements of 
  $(\D,\mathcal{B}(\D))$\,.
Weak convergence of
  probability measures on $\D$ and convergence in distribution of
  stochastic processes  are
 understood with respect to the Skorohod topology.

We recall that a filtered probability space, or a stochastic basis,
$(\Omega,\mathcal{F},\mathbf{F},\mathbf{P})$ is defined as a
probability space $(\Omega,\mathcal{F},\mathbf{P})$ endowed with an increasing
right-continuous   flow $\mathbf{F}=(\mathcal{F}(t),\,t\in\R_+)$ of
   sub-$\sigma$-algebras of $\mathcal{F}$. Such a flow is also
 referred to as a filtration. We will assume without further mention
 that all $\sigma$-algebras we consider are complete with respect to
 the corresponding probability measure.
For  the background on the general theory of stochastic processes, the reader is referred to Jacod and Shiryaev
\cite{JacShi03} and Liptser and Shiryaev \cite{lipshir}\,. For
the background on L\'evy, 
see Bertoin \cite{Ber96}.
For the
properties of regularly and slowly varying functions, see Bingham, Goldie, and
Teugels \cite{Reg89}.



\section{Convergence of the conditioned random walk}
\label{sec:conv-cond-rand}

Given $\alpha\in(0,1)$\,, let $Y=(Y(t)\,, t\in\R_+)$ be the positive
selfsimilar Markov process with the Lamperti representation (see
Lamperti \cite{Lamp72})
\begin{equation}
  \label{eq:1}
Y(t)=\exp\bl(-L(\sigma(t))\br)\,,
\end{equation}
where $L=(L(t)\,,t\in\R_+)$ is a subordinator with
the L\'evy measure $\ind_{\{x>0\}}\alpha e^{-\alpha
  x}(1-e^{-x})^{-\alpha-1}\,dx$\, and $L(0)=0$, and $    \sigma(t)=\inf\{s:\,\int_0^s
\exp\bl(-\alpha L(q)\br)\,dq>t\}\,$.
Processes similar to $Y$ have been investigated in
Chaumont \cite{Cha96} and in Caballero and 
Chaumont \cite{CabCha06} who variously refer to them as 
processes
``conditioned to die at zero'' and
processes ``conditioned to hit zero continuously''.
Although the results of those papers do not apply directly to the
process $Y$ because the authors excluded the case of subordinators,
their methods can be used to study its
  properties. 
We also
mention that the Lamperti representation (\ref{eq:1}) can be
complemented with a representation in terms of the Dol\'eans
exponential:
\begin{equation*}
    Y(t)=\prod_{0<s\le \sigma_1(t)}(1-\Delta L_1(s))\,,
\end{equation*}
where $L_1=(L_1(t)\,,t\in\R_+)$ is a subordinator
  with the L\'evy measure $\ind_{\{x\in(0,1)\}}(1-x)^{\alpha-1}\alpha
  x^{-\alpha-1}\,dx$ and
$    \sigma_1(t)=\inf\{s:\,\int_0^s\prod_{0<s\le q}(1-\Delta
L_1(s))^\alpha\,dq>t\}\,$\,. This is a consequence of the fact that
the process $\bl(-\sum_{0<s\le t} \ln(1-\Delta L_1(s))\,,t\in\R_+\br)$ is distributed as
the process $L$\,.

We define the process $\widetilde{X}=(\widetilde{X}(t)\,,t\in\R_+)$ by letting 
$\widetilde{X}(t)=1-Y(t)$\,. It  is a pure jump
process, with nondecreasing
trajectories, and
 values in $[0,1]$. 
The  predictable jump measure of $\widetilde{X}$ is of the form
$(\nu(\widetilde{X};dt,dx))$\,, where, for 
$G\in\mathcal{B}(\R)$ and $t\in\R_+$,
 \begin{equation}
  \label{eq:73}
  \nu(\bx;[0,t],G)= \int_0^t\int_{G}
\ind_{\{0<x<1-\bx(s)\}}
\Bl(1-\frac{x}{1-\bx(s)}\Br)^{\alpha-1}
\alpha x^{-\alpha-1}
\,dx\,ds\,
\end{equation}
if $\bx\in\D_\uparrow$
and  $\nu(\bx;[0,t],G)=0$ if $\bx\in\D \setminus \D_\uparrow$\,. 
By reversing Lamperti's construction, one can show that the
distribution of $\widetilde{X}$ is uniquely specified by its jump measure,
meaning that there exists a unique  semimartingale 
$\widetilde{X}=(\widetilde{X}(t)\,,t\in\R_+)$ 
defined on a filtered probability space
$(\widetilde{\Omega},\widetilde{\mathcal{F}},
\widetilde{\mathbf{F}},\widetilde{\mathbf{P}})$
with the $\widetilde{\mathbf{F}}$-predictable
jump measure $(\nu(\widetilde{X};dt,dx))$ and such that
 $\widetilde{X}(t)=\sum_{0<s\le t}\Delta \widetilde{X}(s)$\,.
In addition, it can be shown 
that if
$\widetilde{\tau} = \inf \{ t\ge 0 \ : \ \widetilde{X}(t)=1\}$,
then $\widetilde{\tau}<\infty$ a.s. and $\widetilde{X}(\widetilde{\tau}-)=1$.
Moreover, for $n\in\N$,
\begin{equation*}
  \mathbf{\widetilde E}\widetilde{\tau}^n=n!\prod_{k=1}^n
\Bl(\int_0^1
(1-x^{\alpha k })\alpha x^{\alpha-1}(1-x)^{-\alpha-1}\,dx\Br)^{-1}
\end{equation*}
and $  \mathbf{\widetilde E}e^{c\widetilde{\tau}}\le1/(1-c
  \mathbf{\widetilde E}\widetilde{\tau})$ when
  $c<1/\mathbf{\widetilde E}\widetilde{\tau}$\,.

The next theorem, which is analogous to the results available in
Chaumont \cite{Cha96} and in Caballero and 
Chaumont \cite{CabCha06},   confirms the interpretation of $\widetilde{X}$ as ``a
process conditioned not to overshoot''.
Recall that  $X=(X(t)\,,t\in\R_+)$ represents the stable subordinator
  with L\'evy measure $ \alpha
 x^{-\alpha-1}\,dx\,,$ where $\alpha\in(0,1)$\,.
We also  let $\widehat{X}$ represent the process
 $X$ stopped at $\tau$: $\widehat{X}(t)=X(t\wedge\tau)$\,.
 \begin{theorem}
   \label{the:cond_levy}
The conditional laws of $\widehat{X}$ given the events $X(\tau)\le
1+\varepsilon$, considered as distributions on $\D$, weakly converge as
$\varepsilon\downarrow 0$ to the law of $\widetilde{X}$\,.
 \end{theorem}
The proof of this theorem is obtained by arguments modelled on the
proof of Theorem \ref{the:conv_levy} and will not be given here
because of space constraint considerations.
Instead, we outline the intuition behind the form of the jump measure
of $\widetilde{X}$.
 The intensity  of jumps
of size $x$ of $\widetilde{X}$ at a point $\widetilde{X}(t)=u$, where
$u<u+x<1$, can be obtained by ``conditioning'' the intensity of jumps of
$X$ on the event that the overshoot of $X$ over $1$ is not greater
than $\varepsilon$. In other words, it should be given approximately by the product of
the intensity of jumps of $X$ from $u$ to $u+x$, which is $\alpha
x^{-\alpha-1}$, with the probability for $X$ not to exceed level $1$
by more than $\varepsilon$ when starting at $u+x$ over the probability
that
$X$ does not overshoot $1$
by more than $\varepsilon$ starting at $u$.
As follows
by the results of Dynkin \cite{Dyn55} and Rogozin \cite{Rog71},  the probability for the
process $X$ not to overshoot a level $y>0$ by more than $\varepsilon>0$
is asymptotically equivalent to
$\bl(\sin\pi\alpha/(\pi(1-\alpha))\br)y^{\alpha-1}\varepsilon^{1-\alpha}$
as $\varepsilon\to0$\,. Therefore,  the intensity of jumps
 of $\widetilde{X}$ should be $\alpha x^{-\alpha-1}
(1-x/(1-u))^{\alpha-1}$\,. 

The readers who are interested in a complete and self-contained proof of the result of Theorem
2.1, may find it in an earlier version of this paper,
see http://arxiv.org/abs/0712.2637. It also contains the proof of the main result, but
under a slightly stronger version of condition 2.

We now state the main result of the paper. 
Let  
\begin{align}
  \label{eq:22}
    X^{(r)}(t)&=\frac{1}{r}\sum_{i=1}^{\lfloor     t/(1- F (r))\rfloor}\xi_i,\\
\widehat{\tau}^{(r)}&=\inf\{t:\,X^{(r)}(t)\ge 1\}\,.
\end{align}
Note that  $\widehat{\tau}^{(r)} = (1-F(r))\tau^{(r)}$.
We denote by $\widehat{X}^{(r)}=(\widehat{X}^{(r)}(t)\,,t\in\R_+)$
the process $X^{(r)}$ stopped at $\widehat{\tau}^{(r)}$, i.e.,
$\widehat{X}^{(r)}(t)=X^{(r)}(t\wedge\widehat{\tau}^{(r)})$\,.
\begin{theorem}
  \label{the:conv_levy}
Let the following conditions hold:
\begin{enumerate}
\item the right-hand tail of the distribution function $ F $
  is regularly varying
  at infinity with index $-\alpha$, where $\alpha\in(1/2,1)$,
\item there exist $C>0$\,, $B>0$\,, and $\rho\in(0,1)$ such that, for all $z$ large enough and all
  $y\in[\rho ,1-B/z]$,
  \begin{equation*}
    \frac{1- F (yz)}{1- F (z)}\le 1+ C(1-y)\,.
  \end{equation*}
\end{enumerate}
If, in addition,  $F$ is nonlattice,
then, as $r\to\infty$, the conditional distributions of
the $\widehat{X}^{(r)}$ given $\tau^{(r)}<\infty$
weakly converge to the distribution
of  $\widetilde{X}$.
If, instead,  $F$ is  lattice with span $h$,
then, as $n\to\infty$, where $n\in\N$, the conditional distributions of
the $\widehat{X}^{(nh)}$ given $\tau^{(nh)}<\infty$
weakly converge to the distribution
of  $\widetilde{X}$.
\end{theorem}
\begin{remark}
  The restriction on $\alpha$ to be greater than
$1/2$ is due to the fact that a key element in the proof 
is a local renewal theorem  with infinite (or
non-existing) mean which has been established
 for $\alpha \in
(1/2,1)$ only, see  Erickson \cite{Eri70}.
The requirement that $\alpha\in(1/2,1)$ 
is only used in the proof of
  Lemma~\ref{le:denomin}. The rest of the proof applies to any
  distribution $F$ from the domain of attraction of a stable law with
  index $\alpha\in(0,1)$\,.
\end{remark}
\begin{remark}
Under condition 1, the function
 $\ell(x)=x^\alpha\bl(1- F (x)\br)$\,
 is  slowly varying at infinity, i.e.,
 $\lim_{x\to\infty}\ell(yx)/\ell(x)=1$ for all $y>0$\,.
  According to Karamata's theorem,  see Bingham, Goldie, and Teugels \cite{Reg89}
 or Feller \cite{Fel71}, it admits the representation
$\ell(x)=c(x)\exp(\int_1^x \varepsilon(u)/u\,du)$, where $c(x)\to c>0$
 and $\varepsilon(x)\to0$ as $x\to\infty$\,. If $c(x)$ in this
 representation is a constant, or converges to the limit quickly
 enough, then condition 2 of the theorem holds. An example of a
 lattice distribution which meets the requirements of the
 theorem is provided by $F(x)=1-\bl(\lfloor x\rfloor+1\br)^{-\alpha}$
 for $x\ge0$\,.
\end{remark}
\begin{remark}
The distribution  $F$ is lattice if and only if $\widehat{F}$
  is lattice with the same span, where $\widehat{F}$ denotes the
  distribution function of $\xi_1$ under $\mathbf{P}$.
We note that  $ \widehat{F}(dx)=\exp(-\gamma x) F(dx)$\,.
\end{remark}

\begin{remark}
   A sufficient (but not necessary) condition for 
 $F$ to
  have a regularly varying  right-hand tail with index $-\alpha$
 is for the function $e^{\gamma x} \mathbf{P}(\xi_1>x)$ to  be regularly varying
 with index $-\alpha-1$, see
Korshunov \cite{Kor05} for further comments.
Since  the left-hand tail of $F$
 decays exponentially quickly, its right-hand tail  is regularly
 varying with index $-\alpha$ if and only if $F$ belongs to the domain of attraction of
the spectrally positive stable law with index $\alpha$\,, cf.
Gnedenko and Kolmogorov \cite{GneKol68} or Feller \cite{Fel71}.

\end{remark}
The proof of Theorem \ref{the:conv_levy}   will be obtained by an application of the
following  result, which is a special case of Theorem
IX.3.21 on p.546 in Jacod and Shiryaev \cite{JacShi03}\,.
\begin{lemma}
  \label{le:conv_jump}
Consider a sequence $\Xi^{(n)}$ of $\R$-valued 
 semimartingales 
with trajectories of locally bounded variation 
defined on  filtered probability spaces
$(\Omega^{(n)},\mathcal{G}^{(n)},\mathbf{G}^{(n)},\mathbf{P}^{(n)})$\,.
Suppose that the $\Xi^{(n)}$ are pure jump semimartingales
 meaning that $\Xi^{(n)}(t)=\sum_{0<s\le t}\Delta
 \Xi^{(n)}(s)$ for all $t>0$\,and let $\nu^{(n)}(dt,dx)$ denote  their
 predictable
jump measures. Let
an $\R_+$-valued function $K(y; G )$,
where $y\in\R$ and $ G \in\mathcal{B}(\R)$, be
Borel-measurable in $y$ and
 be a $\sigma$-finite measure on $(\R,\mathcal{B}(\R))$
in $ G $  such that $K(y;\{0\})=0$.
Suppose that the following conditions hold:
\begin{enumerate}
\item $\sup_{y\in\R}\int_{\R}\abs{x}\,K(y;dx)<\infty$,
\item for an arbitrary
$\R$-valued continuous function $g(x),\,x\in\R,$
such that $\abs{g(x)}\le M\abs{x},\,x\in\R,$ with some $M>0$,
the function $\int_\R g(x)\,K(y;dx)$ is continuous in $y$,
\item for arbitrary $\delta>0$, $t>0$, and an
$\R$-valued continuous function $g(x),\,x\in\R,$
such that $\abs{g(x)}\le M\abs{x},\,x\in\R,$ with some $M>0$,
\begin{equation*}
  \lim_{n\to\infty}\mathbf{P}^{(n)}
\Bl(\Big| \int_0^t\int_\R g(x)\,\nu^{(n)}(ds,dx)-
\int_0^t\int_\R g(x)\,K(\Xi^{(n)}(s);dx)\,ds \Big| >\delta\Br)=0,
\end{equation*}
\item the $\Xi^{(n)}(0)$ converge in distribution to a random variable $\Xi_0$
  as $n\to\infty$\,,
\item there exists at most one pure jump  semimartingale
 $\Xi=(\Xi(t)\,,t\in\R_+)$ of locally bounded variation
 with an initial condition
  $\Xi_0$  and with predictable jump measure
  $\nu(dt,dx)=K(\Xi(t);dx)\,dt$\,.
\end{enumerate}
Then the $\Xi^{(n)}$ converge in distribution to $\Xi$.
\end{lemma}



\section{Proof of Theorem 2.2.}
\label{sec:proofs}

The proof of Theorem 2.2 follows from Lemmas 3.1--3.5 below. A more detailed
guidance is provided at the end of this Section. 

For the proof, Lemma~\ref{le:conv_jump} will be used with
\begin{equation}
  \label{eq:55}
  K(y; G )= \int_{ G\setminus\{0\}}  \ind_{\{0<x<1-y\}}
\Bl(1-\frac{x}{1-y}\Br)^{\alpha-1}
\alpha x^{-\alpha-1}\,dx\,
\end{equation}
when $y\in(0,1)$ and $K(y; G )=0$ otherwise.
\begin{lemma}
  \label{le:function_K}
The function $K$ satisfies conditions 1 and 2 of Lemma~\ref{le:conv_jump}.
\end{lemma}
\begin{proof}
  For $g(x)$ as in condition 2 of Lemma~\ref{le:conv_jump} and
  for $y\in(0,1)$, we have 
\begin{equation}
  \label{eq:56}
    \int_\R
    g(x)K(y;dx)= (1-y)^{-\alpha}\int_0^1g(u(1-y))(1-u)^{\alpha-1}
\alpha u^{-\alpha-1}\,du\,.
\end{equation}
Thus, $  \int_\R \abs{x}K(y;dx)\le\alpha
\int_0^1u^{\alpha-1}(1-u)^{-\alpha}\,du,$
so condition 1 of Lemma~\ref{le:conv_jump} holds.
Condition 2   follows by the assumption that
$\abs{g(x)}\le M\abs{x}$, continuity of $g(x)$, and Lebesgue's bounded convergence theorem.
\end{proof}

In order to apply Lemma~\ref{le:conv_jump}, we need to compute the
predictable triplets of the conditioned processes $\hat{X}^{(r)}$\,.
We start with introducing a change of measure. Let the absolutely continuous
with respect to $\bP^\ast$ probabilty measure $\bQ^{(r)}$ on
$(\Omega,\mathcal{F}_{\tau^{(r)}})$  be defined by
\begin{equation}
  \label{eq:3}
  \frac{d\bQ^{(r)}}{d\bP^\ast}=
\frac{e^{-\gamma \chi^{(r)}}}{\bE^\ast
e^{-\gamma \chi^{(r)}}}\,.
\end{equation}
By \eqref{eq:66},
 for $\Gamma\in \mathcal{F}_{\tau^{(r)}}$,
\begin{equation*}
    \bP(\Gamma|\tau^{(r)}<\infty)=
\frac{\bE^\ast e^{-\gamma \chi^{(r)}}\ind_\Gamma}{\bE^\ast e^{-\gamma \chi^{(r)}}}\,,
\end{equation*}
so,  for a Borel subset $U$ of $\D$,
\begin{equation}
  \label{eq:6}
\bP(\widehat{X}^{(r)}\in
 U|\tau^{(r)}<\infty)=  \bQ^{(r)}(\widehat{X}^{(r)}\in U).
\end{equation}
Let $\mu^{(r)}$ denote the jump measure of $X^{(r)}$,  so that
\begin{equation}
  \label{eq:10}
  \mu^{(r)}([0,t], G )=\sum_{i=1}^{\lfloor  t/(1- F (r))\rfloor}
\ind_{\{\xi_i/r\in  G \setminus\{0\}\}},
\end{equation}
where $ G $ is a Borel subset of $\R$. Clearly,
\begin{equation}
  \label{eq:11}
  X^{(r)}(t)=\int_0^t\int_\R x\,\mu^{(r)}(ds,dx)\,.
\end{equation}
Let  $\mathcal{F}^{(r)}(t)$ represent  the $\sigma$-algebra generated by the random variables
$X^{(r)}(s)$ for $s\le t\,.$ 
Let $\mathbf{F}^{(r)}=(\mathcal{F}^{(r)}(t)\,,t\in\R_+)$
be the filtration associated with $X^{(r)}$ and
let $\widehat{\mathbf{F}}^{(r)}=(\widehat{\mathcal{F}}^{(r)}(t)\,,t\in\R_+)$,
with
$\widehat{\mathcal{F}}^{(r)}(t)=\mathcal{F}^{(r)}(t\wedge\widehat{\tau}^{(r)})$,
be the filtration associated with $\widehat{X}^{(r)}$\,.
The $\mathbf{F}^{(r)}$-predictable
jump measure of $X^{(r)}$  under  $\bP^\ast$ is of the form
\begin{equation}
  \label{eq:9}
  \nu^{(r)}([0,t], G )= \bigg\lfloor\, \frac{ t}{1- F (r)}\,\bigg \rfloor
\int_{ G \setminus\{0\}}  F (r\,dx)\,.
\end{equation}
  The jump measure
 of $\widehat{X}^{(r)}$ is given by
\begin{equation}
  \label{eq:18}
  \widehat{\mu}^{(r)}([0,t], G )=\sum_{i=1}^{\lfloor  (t\wedge
    \widehat{\tau}^{(r)})/(1- F (r))\rfloor}
\ind_{\{\xi_i/r\in  G \setminus\{0\}\}},
\end{equation}
and the $\widehat{\mathbf{F}}^{(r)}$-predictable jump measure 
 of $\widehat{X}^{(r)}$ under $\mathbf{P}^\ast$ is 
\begin{equation}
  \label{eq:25}
\widehat{\nu}^{(r)}([0,t], G )=
 \nu^{(r)} ([0,t\wedge     \widehat{\tau}^{(r)}],G)\,.
\end{equation}
For $y>0$, we let
\begin{equation}
  \label{eq:35}
    u^{(r)}(y)=\bE^\ast e^{-\gamma \chi^{(ry)}}\,.
\end{equation}
\begin{lemma}
  \label{le:pred-meas}
There exists a version of the $\widehat{\mathbf{F}}^{(r)}$-predictable
jump measure of $\widehat{X}^{(r)}$  under  $\mathbf{Q}^{(r)}$
 of the form
\begin{multline*}
  \widetilde{\nu}^{(r)}(dt,dx)=\ind_{\{\widehat{X}^{(r)}(t-)<1\}}\Bl[
\ind_{\{\widehat{X}^{(r)}(t-)+x\ge 1\}}\,\frac{e^{-\gamma
      r(\widehat{X}^{(r)}(t-)+x- 1)}}{u^{(r)}(1-\widehat{X}^{(r)}(t-))}
\\+\ind_{\{\widehat{X}^{(r)}(t-)+x<1\}}\,\frac{u^{(r)}
(1-\widehat{X}^{(r)}(t-)-x)}{u^{(r)}(1-\widehat{X}^{(r)}(t-))}\Br]\,\nu^{(r)}(dt,dx)\,.
\end{multline*}
\end{lemma}
\begin{proof}
Let us introduce
\begin{equation}
  \kappa^{(r)}(y)=\inf\{t\in\R_+:\, X^{(r)}(t)\ge y\}
\end{equation}
and
\begin{equation}
  \label{eq:5}
  \eta^{(r)}(y)=X^{(r)}(\kappa^{(r)}(y))-y\,.
\end{equation}
Note that $\widehat{\tau}^{(r)}=\kappa^{(r)}(1)$ and
 $\chi^{(ry)}=ry\eta^{(ry)}(1)=r\eta^{(r)}(y)$\,,
so by   \eqref{eq:35}
\begin{equation}
  \label{eq:13}
  u^{(r)}(y)=\bE^\ast e^{-\gamma    r\eta^{(r)}(y)}\,.
\end{equation}
Also, by \eqref{eq:3} the density process of $\bQ^{(r)}$ with respect
to $\bP^\ast$, which we denote by $Z^{(r)}=(Z^{(r)}(t)\,,t\in\R_+)$, can be written as
\begin{equation*}
  Z^{(r)}(t)=\bE^\ast\Bl[\frac{e^{-\gamma r\eta^{(r)}(1)}}{\bE^\ast
e^{-\gamma r\eta^{(r)}(1)}}\big|\, \mathcal{F}^{(r)}(t\wedge \kappa^{(r)}(1) )\Br]\,.
\end{equation*}

Since $\eta^{(r)}(1)$ is
an $\mathcal{F}^{(r)}(\kappa^{(r)}(1))$-measurable random variable
and  $X^{(r)}$ has independent stationary increments,
\begin{multline*}
  Z^{(r)}(t)=\ind_{\{\kappa^{(r)}(1)\le t\}}
\bE^\ast\Bl[\frac{e^{-\gamma r\eta^{(r)}(1)}}{\bE^\ast
e^{-\gamma r\eta^{(r)}(1)}}\big|\, \mathcal{F}^{(r)}(\kappa^{(r)}(1))\Br]+
\ind_{\{\kappa^{(r)}(1)> t\}}
\bE^\ast\Bl[\frac{e^{-\gamma r\eta^{(r)}(1)}}{\bE^\ast
e^{-\gamma r\eta^{(r)}(1)}}\big|\, \mathcal{F}^{(r)}(t)\Br]\\=
\ind_{\{\kappa^{(r)}(1)\le t\}}\frac{e^{-\gamma r\eta^{(r)}(1)}}{\bE^\ast
e^{-\gamma r\eta^{(r)}(1)}}+
\ind_{\{\kappa^{(r)}(1)> t\}}\frac{\bE^\ast e^{-\gamma
    r\eta^{(r)}(y)}\big|_{y=1-X^{(r)}(t)}}{
\bE^\ast e^{-\gamma r\eta^{(r)}(1)}}\,.
\end{multline*}
Thus, on recalling \eqref{eq:13} and taking into account that
$\widehat{X}^{(r)}(t)=X^{(r)}(t)$ for $t<\kappa^{(r)}(1)$,
\begin{equation}
  \label{eq:16}
  Z^{(r)}(t)=\ind_{\{\kappa^{(r)}(1)\le t\}}
\frac{e^{-\gamma r\eta^{(r)}(1)}}{u^{(r)}(1)}+
\ind_{\{\kappa^{(r)}(1)> t\}}\frac{u^{(r)}(1-\widehat{X}^{(r)}(t))}{u^{(r)}(1)}\,.
\end{equation}
By Liptser and Shiryaev \cite[p.223]{lipshir} or Jacod and Shiryaev
\cite[p.170]{JacShi03},  the $\widehat{\mathbf{F}}^{(r)}$-predictable
jump measure of $\widehat{X}^{(r)}$ under
$\bQ^{(r)}$ admits a version of the form
\begin{equation}
  \label{eq:15}
  \widetilde{\nu}^{(r)}(dt,dx)=Y^{(r)}(t,x)\widehat{\nu}^{(r)}(dt,dx),
\end{equation}
where
\begin{equation}
  \label{eq:17}
  Y^{(r)}(t,x)=\frac{\ind_{\{Z^{(r)}(t-)>0\}}}{Z^{(r)}(t-)}\,\mathbf{M}^{\mathbf{P}^\ast}_{\widehat{\mu}^{(r)}}
(Z^{(r)}|\widetilde{\mathcal{P}})(t,x).
\end{equation}
For the reader's convenience, we recall  that $\widetilde{\mathcal{P}}$ denotes the $\sigma$-algebra on
$\Omega\times\R_+\times\R$ that is the product of the predictable
$\sigma$-algebra on $\Omega\times\R_+$ associated with
$\widehat{\mathbf{F}}^{(r)}$
 and the Borel $\sigma$-algebra
on $\R$, and $\mathbf{M}^{\mathbf{P}^\ast}_{\widehat{\mu}^{(r)}}$ denotes the
measure on $\Omega\times\R_+\times\R
$ defined by the equality
\begin{equation*}
  \mathbf{M}^{\mathbf{P}^\ast}_{\widehat{\mu}^{(r)}}f=
\mathbf{E}^\ast\int_0^\infty\int_\R f(\omega,t,x)\,\widehat{\mu}^{(r)}(dt,dx)\,
\end{equation*}
for $f\ge0$\,, where $\widehat{ \mu}^{(r)}$ is the  jump measure
 of $\widehat{ X}^{(r)}$,
i.e.,
\begin{equation}
  \label{eq:82}
  \widehat{\mu}([0,t], G )=\sum_{0<s\le t\wedge\tau}\ind_{\{\Delta X(s)\in G\setminus\{0\} \}}\,.
\end{equation}
Accordingly, $\mathbf{M}^{\mathbf{P}^\ast}_{\widehat{\mu}^{(r)}}
(Z^{(r)}|\widetilde{\mathcal{P}})(t,x)$ is the conditional expectation of
$Z^{(r)}$ with respect to $\widetilde{\mathcal{P}}$, i.e., it is a
$\widetilde{\mathcal{P}}$-measurable function $g(\omega,t,x)$ such that
$  \mathbf{M}^{\mathbf{P}\ast}_{\widehat{\mu}^{(r)}}hZ^{(r)}=
  \mathbf{M}^{\mathbf{P}^\ast}_{\widehat{\mu}^{(r)}}hg$\, for all
  nonnegative $\widetilde{\mathcal{P}}$-measurable functions $h$\,.

In order to find $\mathbf{M}^{\mathbf{P}^\ast}_{\widehat{\mu}^{(r)}}
(Z^{(r)}|\widetilde{\mathcal{P}})(t,x)$, we write by \eqref{eq:18}, \eqref{eq:16}, and the
definition of $\eta^{(r)}(1)$ in \eqref{eq:5}, for a
$\widetilde{\mathcal{P}}$-measurable nonnegative function $h$,
\begin{multline*}
\mathbf{M}^{\mathbf{P}^\ast}_{\widehat{\mu}^{(r)}}h Z^{(r)}=
\mathbf{E}^\ast  \int_0^\infty\int_\R h(t,x)\,
\ind_{\{\widehat{X}^{(r)}(t-)<1\}}\Bl(\ind_{\{\widehat{X}^{(r)}(t-)+x\ge1\}}
\frac{e^{-\gamma r(\widehat{X}^{(r)}(t-)+x- 1)}}{u^{(r)}(1)}\\+
\ind_{\{\widehat{X}^{(r)}(t-)+x<1\}}\frac{u^{(r)}(1-\widehat{X}^{(r)}(t-)-x)}{u^{(r)}(1)}\Br)
\,\widehat{\mu}^{(r)}(dt,dx)\,.
\end{multline*}
The expression in parentheses is $\widetilde{\mathcal{P}}$-measurable.
Hence,
\begin{multline*}
  \mathbf{M}^{\mathbf{P}^\ast}_{\widehat{\mu}^{(r)}}
(Z^{(r)}|\widetilde{\mathcal{P}})(t,x)=
\frac{\ind_{\{\widehat{X}^{(r)}(t-)<1\}}}{u^{(r)}(1)}\,\Bl(\ind_{\{\widehat{X}^{(r)}(t-)+x\ge1\}}
e^{-\gamma r(\widehat{X}^{(r)}(t-)+x- 1)}\\+
\ind_{\{\widehat{X}^{(r)}(t-)+x<1\}}u^{(r)}(1-\widehat{X}^{(r)}(t-)-x)\Br)\,.
\end{multline*}
The expression for $\widetilde{\nu}^{(r)}(dt,dx)$ in the statement of the
lemma follows now  by \eqref{eq:25}, \eqref{eq:16}, \eqref{eq:15}, and \eqref{eq:17}.
\end{proof}
In what follows, we work with the version of the  $\widehat{\mathbf{F}}^{(r)}$-predictable
jump measure of $\widehat{X}^{(r)}$  given
 in the statement of Lemma~\ref{le:pred-meas}.
We  study the properties of $\widetilde{\nu}^{(r)}$\,.
Let
\begin{subequations}
  \begin{align}
    \label{eq:34b}
    \widetilde{\nu}^{(r)}_1(dt,dx)&=\ind_{\{\widehat{X}^{(r)}(t-)<1\}}
\ind_{\{\widehat{X}^{(r)}(t-)+x\ge 1\}}\,\frac{e^{-\gamma
      r(\widehat{X}^{(r)}(t-)+x- 1)}}{u^{(r)}(1-\widehat{X}^{(r)}(t-))}\,\nu^{(r)}(dt,dx)
\,,\\    \label{eq:34a}
    \widetilde{\nu}^{(r)}_2(dt,dx)&=
\ind_{\{\widehat{X}^{(r)}(t-)<1\}}\ind_{\{x<0\}}\,
\frac{u^{(r)}
(1-\widehat{X}^{(r)}(t-)-x)}{u^{(r)}(1-\widehat{X}^{(r)}(t-))}
\,\nu^{(r)}(dt,dx)\,,\\
    \label{eq:34c}
    \widetilde{\nu}^{(r)}_3(dt,dx)&=\ind_{\{x>0\}}\,
\ind_{\{\widehat{X}^{(r)}(t-)+x<1\}}\,
\frac{u^{(r)}
(1-\widehat{X}^{(r)}(t-)-x)}{u^{(r)}(1-\widehat{X}^{(r)}(t-))}
\,\nu^{(r)}(dt,dx)\,.
  \end{align}
\end{subequations}
Note that
\begin{equation}
  \label{eq:29}
  \widetilde{\nu}^{(r)}=\widetilde{\nu}^{(r)}_1+\widetilde{\nu}^{(r)}_2+\widetilde{\nu}^{(r)}_3\,.
\end{equation}
One can see that $\widetilde{\nu}^{(r)}_1$ describes the intensity
of  jumps of $\widehat{X}^{(r)}$ reaching level 1, that 
$\widetilde{\nu}^{(r)}_2$  concerns  the intensity
of downward jumps, and that 
$\widetilde{\nu}^{(r)}_3$ characterises the intensity
of upward jumps.
Not unexpectedly, the first two measures are inconsequential, as the
next lemma shows.
\begin{lemma}
\label{le:conv}
Let the hypotheses  of Theorem~\ref{the:conv_levy}  hold. Then,
 for $i=1,2$,
\begin{equation*}
\lim_{r\to\infty}\sup_{\omega\in\Omega}  \int_0^{t}\int_\R
\abs{x}\,\widetilde{\nu}_i^{(r)}(ds,dx)
=0\,
\end{equation*}
 if $F$ is  nonlattice  and
\begin{equation*}\lim_{n\to\infty}
\sup_{\omega\in\Omega}  \int_0^{t}\int_\R
\abs{x}\,\widetilde{\nu}_i^{(nh)}(ds,dx)
=0\,
\end{equation*}
if $F$ is  lattice  with span $h$\,.
\end{lemma}
\begin{proof}
Suppose that $F$ is  nonlattice.
We start with $i=1$.  Using \eqref{eq:9} and integrating on $x$
yields on taking into account that $\widehat{X}^{(r)}(s)$ is constant on the intervals
$[i(1- F (r)),(i+1)(1- F (r)))$,
\begin{multline}
  \label{eq:24}
\int_0^{t}\int_\R
\abs{x}\,\widetilde{\nu}_1^{(r)}(ds,dx)\le
    \int_0^{t}\int_\R
\ind_{\{\widehat{X}^{(r)}(s-)<1,\,\widehat{X}^{(r)}(s-)+x\ge 1\}}\,
\frac{e^{-\gamma
      r(\widehat{X}^{(r)}(s-)+x- 1)}}{u^{(r)}
(1-\widehat{X}^{(r)}(s-))}\,\abs{x}\,\nu^{(r)}(ds,dx)
\\
=\int_0^{t} \frac{\ind_{\{\widehat{X}^{(r)}(s-)<1\}}}{u^{(r)}
(1-\widehat{X}^{(r)}(s-))}
\int_{1-\widehat{X}^{(r)}(s-)}^\infty x\,
e^{-\gamma       r(\widehat{X}^{(r)}(s-)+x- 1)}
\, F (r\,dx)\,
d\, \lfloor (1- F (r))^{-1}\, s \rfloor
\\
=\frac{1}{1- F (r)}\,\int_0^{\lfloor  t/(1- F (r))\rfloor(1- F (r))}
\frac{\ind_{\{\widehat{X}^{(r)}(s)<1\}}}{u^{(r)}
(1-\widehat{X}^{(r)}(s))}\int_{1-\widehat{X}^{(r)}(s)}^\infty x\,
e^{-\gamma       r(\widehat{X}^{(r)}(s)+x- 1)}
\, F (r\,dx)\,d s\\
\le\frac{t}{1- F (r)}\,\sup_{y\in[0,1]}
\frac{1}{u^{(r)}(y)}\int_y^\infty x\,
e^{-\gamma       r(x-y)}
\, F (r\,dx)\,.
\end{multline}
For $y\in[ 0,1]$ and $A\in(0,r)$, on recalling \eqref{eq:13},
 employing a change of variables and
integrating by parts,
\begin{multline}
  \label{eq:26}
    \frac{1}{(1- F (r))u^{(r)}(y)}\int_{y}^\infty x\,
e^{-\gamma       r(x- y)}\, F (r\,dx)=
  \frac{1}{r(1- F (r))\bE^\ast e^{-\gamma \chi^{(ry)}}}\int_{ry}^\infty x\,
e^{-\gamma       (x- ry)}\, F (dx)\\
\le \ind_{\{ry\le A\}}
\, \frac{1}{r(1- F (r))\displaystyle\inf_{z\le A}\bE^\ast e^{-\gamma \chi^{(z)}}}\,
\sup_{z\le A}\int_{A}^\infty x\,
e^{-\gamma       (x- z)}\, F (dx)
\\+\ind_{\{ry> A\}}\,\frac{1}{ry(1- F (ry))
\bE^\ast e^{-\gamma \chi^{(ry)}}}\,
\frac{y(1- F (ry))}{1- F (r)}\,
\gamma\int_0^\infty e^{-\gamma z}\int_{ry}^{z+ry} x\,
\, F (dx)\,dz.
\end{multline}
We first work with the second term on the rightmost side.
By the first assertion of part 2 of Lemma~\ref{le:slowly}, given arbitrary $\varepsilon>0$,
$ \sup_{y\in[A/r,1]}y(1- F (ry))/(1- F (r))\le1+\varepsilon$ for
all $A$ and $r$ large enough. By Lemma~\ref{le:denomin}, $(ry)(1- F (ry))
\bE^\ast \exp(-\gamma \chi^{(ry)})$ converges to a positive limit as
$ry\to\infty$\,.

Next,  for $u>0$,
\begin{equation}
  \label{eq:77}
    \int_{u}^{z+u} x\,  F (dx)\le
  (z+u)\bl( F (z+u)- F (u)\br)\,.
\end{equation}
Condition 2 of Theorem~\ref{the:conv_levy} implies that,
 for all $a>B$,
  $ \lim_{x\to\infty}x( F (x+a)- F (x))=0$\,,
so
\begin{equation*}
  \lim_{u\to\infty}(z+u)\bl( F (z+u)- F (u)\br)=0\,.
\end{equation*}
Also, for $u\ge1$,
\begin{multline*}
  e^{-\gamma z}(z+u)\bl( F (z+u)- F (u)\br)\le
 e^{-\gamma z}z+
 e^{-\gamma z}u\bl( F (z+u)- F (u)\br)
\ind_{\{z\le 2\ln u/\gamma\}}+
  e^{-\gamma z}u\ind_{\{z>2 \ln u/\gamma\}}\\
\le e^{-\gamma z}z+e^{-\gamma z}u\bl( F \bl(
\,\frac{2\ln u}{\gamma}+u\br)- F (u)\br)
+  e^{-\gamma z/2}\,.
\end{multline*}
Another application of condition 2 shows that,
 for all $a>0$, $ \lim_{x\to\infty}
x( F (x+a\ln x)- F (x))=0$\,.
Thus, the rightmost side of the latter display
is bounded from above by
a function of $z$ with a finite   $(0,\infty)$-integral  with respect
 to Lebesgue measure\,. Hence, by
the dominated convergence theorem,
\begin{equation*}
  \lim_{u\to\infty}\gamma\int_0^\infty e^{-\gamma z}\int_{u}^{u+z} x\,
\, F (dx)\,dz=0.
\end{equation*}
Thus,  given arbitrary $\varepsilon>0$,
 if $A$ is large enough, then the
second term on the rightmost side of \eqref{eq:26} is less than $\varepsilon$.
Since the first term on the rightmost side  of \eqref{eq:26} converges to zero
as $r\to\infty$ for  fixed $A$, we conclude that
\begin{equation*}
  \lim_{r\to\infty}\sup_{y\in[0,1]}\frac{1}{(1- F (r))u^{(r)}(y)}\int_{y}^\infty x\,
e^{-\gamma       r(x- y)}\, F (r\,dx)=0\,.
\end{equation*}
By \eqref{eq:24},
$\lim_{r\to\infty}\sup_{\omega\in\Omega}  \int_0^{t}\int_\R
\abs{x}\,\widetilde{\nu}_1^{(r)}(ds,dx)=0$\,.

Let $i=2$. By analogy   with  \eqref{eq:24} and \eqref{eq:26}
  and by
 the bound $u^{(r)}(y)\le1$,
  we have that for $A>0$,
\begin{multline*}
  \int_0^t\int_\R \abs{x}\,\widetilde{\nu}^{(r)}_2(ds,dx)
\le\frac{1}{1- F (r)}\int_0^t
\int_\R\ind_{\{\widehat{X}^{(r)}(s)<1\}}\,x^-
\,\frac{u^{(r)}
(1-\widehat{X}^{(r)}(s)-x)}{u^{(r)}(1-\widehat{X}^{(r)}(s))}
\,\, F (r\,dx)\,d s\\
\le
\,\frac{1}{r(1- F (r))\displaystyle\inf_{z\le A}\bE^\ast e^{-\gamma \chi^{(z)}}}
\int_0^t
\int_\R\ind_{\{\widehat{X}^{(r)}(s)<1\}}\ind_{\{r(1-\widehat{X}^{(r)}(s))\le A\}}\,x^-
\,\, F (dx)\,d s\\
+\frac{1}{r(1- F (r))}\,
\frac{\sup_{z\ge A}\mathbf{E}^\ast \exp(-\gamma\chi^{(z)})}{
\inf_{z\ge A}\mathbf{E}^\ast \exp(-\gamma\chi^{(z)})}\int_0^t
\int_\R\ind_{\{\widehat{X}^{(r)}(s)<1\}}\ind_{\{r(1-\widehat{X}^{(r)}(s))> A\}}\,x^-
\, F (dx)\,ds\,.
\end{multline*}
Thus, by Lemma~\ref{le:denomin}, given $\varepsilon>0$, for all
$A$ large enough and all $r>A$,
\begin{equation*}
  \int_0^t\int_\R \abs{x}\,\widetilde{\nu}^{(r)}_2(ds,dx)
\le
\frac{t}{r(1- F (r))}\,
\Bl( \frac{1}{\displaystyle\inf_{z\le A}\bE^\ast e^{-\gamma \chi^{(z)}}}+1+\varepsilon\Br)
\,\int_\R\,x^-\, F (dx)\,.
\end{equation*}

Since $\int_\R\,x^-\, F (dx)=\int_{-\infty}^0 (-x)\exp(\gamma
x)\,\widehat{F}(dx)<\infty$,
the function $r^\alpha(1- F (r))$ is slowly varying at infinity, and
  $\alpha<1$, we conclude that
$\sup_{\omega\in\Omega}  \int_0^t\int_\R \abs{x}\,\widetilde{\nu}^{(r)}_1(ds,dx)$ converges
 to $0$ as $r\to\infty$.

The assertion of
 the lemma for nonlattice $F$  has been proved.
The proof for the case where $F$ is lattice 
with span $h$ proceeds analogously. The only changes consist in
assuming that $r=nh$, that $y$ is of the form $k/n$ for
$k=1,2,\ldots,n$, that $z$ is of the form $kh$ for $k=1,2,\ldots$, and
 that $n\to\infty$,
 and in using the part of Lemma~\ref{le:denomin} that concerns the lattice case.
\end{proof}
We are now in a position to verify condition 3 of Lemma~\ref{le:conv_jump}.
\begin{lemma}
  \label{le:conv_char}
Let the hypotheses of Theorem~\ref{the:conv_levy} hold. 
Let $K$ be defined by \eqref{eq:55},  and  
$g(x)$ be a  bounded continuous
function  such that $\abs{g(x)}\le M\abs{x}$ for some $M>0$.
Take any $t>0$ and $\varepsilon>0$.\\ 
If $F$ is  nonlattice, then
\begin{equation*}
  \lim_{r\to\infty}
\mathbf{Q}^{(r)}\Bl(\Big|\int_0^t\int_\R g(x)\,\widetilde{\nu}^{(r)}(ds,dx)-
\int_0^t\int_\R g(x)\,K(\widehat{X}^{(r)}(s);dx)\,ds\Big|>\varepsilon\Br)=0\,.
\end{equation*}
If, instead, $F$ is  lattice  with span $h$, then
\begin{equation*}
  \lim_{n\to\infty}
\mathbf{Q}^{(nh)}\Bl(\Big|\int_0^t\int_\R g(x)\,\widetilde{\nu}^{(nh)}(ds,dx)-
\int_0^t\int_\R g(x)\,K(\widehat{X}^{(nh)}(s);dx)\,ds\Big|>\varepsilon\Br)=0\,.
\end{equation*}
\end{lemma}
\begin{proof}
Suppose $F$ is nonlattice.
We will prove that $\lim_{r\to\infty} \sup_{\omega\in\Omega} M_2=0$ where
\begin{equation}
  \label{eq:28}
 M_2 =
\Big|\int_0^t\int_\R g(x)\,\widetilde{\nu}^{(r)}(ds,dx)-
\int_0^t\int_\R g(x)\,K(\widehat{X}^{(r)}(s);dx)\,ds\Big|\,.
\end{equation}
By Lemma~\ref{le:conv},
\begin{equation}
  \label{eq:33}
\lim_{r\to\infty}  \sup_{\omega\in\Omega}\int_0^t\int_\R
\abs{g(x)}\,\widetilde{\nu}^{(r)}_i(ds,dx)=0\,,i=1,2\,.
\end{equation}
We  turn our attention to $\widetilde{\nu}^{(r)}_3$.
By \eqref{eq:34c},
 \begin{equation*}
\widetilde{\nu}^{(r)}_3(dt,dx)=
\ind_{\{0<x<1-\widehat{X}^{(r)}(t-)\}}\,\frac{u^{(r)}
(1-\widehat{X}^{(r)}(t-)-x)}{u^{(r)}(1-\widehat{X}^{(r)}(t-))}
\,\nu^{(r)}(dt,dx) \,.
\end{equation*}
Let $h(t,r)= \lfloor  t/(1- F (r))\rfloor(1- F (r))$.
By \eqref{eq:9},
$$
\int_0^t\int_\R g(x)\,\widetilde{\nu}^{(r)}_3(ds,dx)  =
\int_0^{h(t,r)}\int_0^1 V_1^{(r)}(\widehat{X}^{(r)}(s),x)\,
\,\frac{ F \bl(r(1-\widehat{X}^{(r)}(s))\,dx\br)}{
1- F (r)} \ ds
$$
where
$$
V_1^{(r)}(z,x)=g\bl((1-z)x\br)\ind_{\{z<1\}}\,
\,\frac{u^{(r)}\bl((1-z)(1-x)\br)}{u^{(r)}
(1-z)}.
$$
Let also
\begin{align*}
V_2(z,x)&= g((1-z)x)\ind_{\{z<1\}}\,
(1-z)^{-\alpha}(1-x)^{\alpha -1}ax^{-\alpha -1}
\intertext{and} 
V_3(z,x)&=
g((1-z)x)\ind_{\{z<1\}}\,
(1-x)^{\alpha-1}
.
\end{align*}
Then, for $\eta\in(0,1)$, recalling \eqref{eq:56},
\begin{eqnarray*}
  \label{eq:30}
  M_2 &\leq &
 \int_0^t \int_{1-\eta}^1 \abs{V_2(\widehat{X}^{(r)}(s),x)} \,d x\,d s
+\int_0^{h(t,r)} \int_{1-\eta}^1
\abs{V_1^{(r)}(\widehat{X}^{(r)}(s), x)}\,\frac{ F \bl(r(1-\widehat{X}^{(r)}(s))\,dx\br)}{1- F (r)}\,d
s\, \\
&+&\int_{h(t,r)}^t \int_0^{1-\eta}
\abs{V_3(\widehat{X}^{(r)}(s), x)}\,\frac{ F \bl(r(1-\widehat{X}^{(r)}(s))\,dx\br)}{1- F (r)}\,
d s\\ &
+&\int_0^{h(t,x)} \int_0^{1-\eta}
\abs{V_1^{(r)}(\widehat{X}^{(r)}(s)
,  x)-V_3(\widehat{X}^{(r)}(s),  x)}\,\frac{ F \bl(r(1-\widehat{X}^{(r)}(s))\,dx\br)}{1-
F (r)}
\,d s \\
&+&\Big|\int_0^t \int_0^{1-\eta}
V_3(\widehat{X}^{(r)}(s),  x)\,\frac{ F \bl(r(1-\widehat{X}^{(r)}(s))\,dx\br)}{1- F (r)}\,d s
-
 \int_0^t \int_0^{1-\eta}
V_2(\widehat{X}^{(r)}(s),x)\,dx\,d s\Big|\,.
\end{eqnarray*}
We denote the terms on the right-hand side 
 by $I_1$, $I_2$, $I_3$, $I_4$, and $I_5$, respectively. We treat them
 successively.
We have
 \begin{equation*}
   I_1\le
   M\int_0^t(1-\widehat{X}^{(r)}(s))^{1-\alpha}\ind_{\{\widehat{X}^{(r)}(s)<1\}}\,
\int_{1-\eta}^1(1-x)^{\alpha-1}\,\alpha x^{-\alpha}\,dx\,ds
\le Mt
\int_{1-\eta}^1(1-x)^{\alpha-1}\,\alpha x^{-\alpha}\,dx\,,
 \end{equation*}
so
\begin{equation}
  \label{eq:32}
  \lim_{\eta\to0}\limsup_{r\to\infty}\sup_{\omega\in\Omega}I_1=0\,.
\end{equation}

We now prove that
\begin{equation}
  \label{eq:31}
  \lim_{\eta\to0}\limsup_{r\to\infty}\sup_{\omega\in\Omega}I_2=0\,.
\end{equation}
By Lemma \ref{le:denomin},
$ry(1- F (ry)) u^{(r)}(y)\to C_0>0$ as $ry\to\infty$  and
$\sup_{ry\ge0} ry(1- F (ry)) u^{(r)}(y)<\infty$\,. On denoting the latter
supremum  by $N$, we have that,  for
 $r y $ is large enough and for $x\in(0,1)$,
 \begin{equation}
   \label{eq:2}
   \frac{u^{(r)}\bl( y (1-x)\br)}{u^{(r)}
 (y) }  \le \frac{2N}{C_0} \,\frac{1- F (r
 y )}{(1-x)(1- F (r
 y (1-x)))}\,.
\end{equation}
Since  $u^{(r)}\bl( y (1-x)\br)\le1$, we can
write, for large enough $A\in(0,r)$ and for $A'\in(B,A)$,
\begin{multline}
  \label{eq:39}
    I_2 \le  M\int_0^t \int_{1-\eta}^1
\ind_{\{\widehat{X}^{(r)}(s)<1\}}x\,\bl(1-\widehat{X}^{(r)}(s))\br)
\,\frac{u^{(r)}\bl((1-\widehat{X}^{(r)}(s))(1-x)\br)}{u^{(r)}
(1-\widehat{X}^{(r)}(s))}\,
\frac{ F \bl(r(1-\widehat{X}^{(r)}(s))\,dx\br)}{1- F (r)}\,ds
\\\le  Mt\sup_{y\in[0,1]} \int_{1-\eta}^1
xy\,\,\frac{u^{(r)}\bl(y(1-x)\br)}{u^{(r)}(y)}\,
\frac{ F \bl(ry\,dx\br)}{1- F (r)}\,ds \quad \ \ \ \quad \quad \quad \quad
\quad \quad \quad \quad \quad \quad \quad \quad \quad \quad \quad \quad
\\\le  Mt\,\frac{A}{r(1- F (r))}\,
\frac{ F (A)}{\displaystyle\inf_{z\le A}\mathbf{E}^\ast \exp(-\gamma\chi^{(z)})}\,
+
Mt\sup_{y\in(A/r,1]} \int_{1-\eta}^1
\ind_{\{ry(1-x)\le A'\}}
\,\frac{x\,y}{u^{(r)}(y)}\,
\frac{ F (ry\,dx)}{1- F (r)}\,
\\+
Mt\,\frac{2N}{C_0}\,
\sup_{y\in(A/r,1]}\int_{1-\eta}^1
\ind_{\{ry(1-x)>A'\}}\,xy\,\frac{1- F (ry)}{(1-x)(1- F (ry(1-x)))}
\frac{ F (ry\,dx)}{1- F (r)}\,. \quad \quad \quad \quad
\end{multline}
The first term on the rightmost side of \eqref{eq:39}  tends to zero as $r\to\infty$\,.
By the fact that $1-x\le A'/(ry)$ when the integrand of the second term is
positive and by \eqref{eq:2}, that term is not greater than
\begin{equation*}
  Mt\,\frac{2}{C_0}\,
\sup_{y\in(A/r,1]}\frac{y(1- F (ry))}{1- F (r)}
\sup_{v>A}v( F (v)- F (v-A'))\,.
\end{equation*}
By the first assertion of  part 2
 of Lemma~\ref{le:slowly}, the supremum in the middle  is bounded from above for all
large enough $A$ and $r$\,. By condition 2 of Theorem~\ref{the:conv_levy},
the other  supremum tends to zero as $A\to\infty$\,.
Thus, the second term on the rightmost side of \eqref{eq:39} tends to
 zero as $A\to\infty$ and $r\to\infty$\,.

Let us consider  the third term on the rightmost side of \eqref{eq:39}.
Pick arbitrary $\varepsilon\in(0,\alpha\wedge(1-\alpha))$\,.
The function $\ell(x)=x^{\alpha}(1- F (x))$ is slowly varying at
infinity, so by part 1 of Lemma~\ref{le:slowly}, for
  $A'$ large enough and for $ry(1-x)>A'$,
\begin{equation*}
  \frac{1- F (ry)}{1- F (ry(1-x))}=
(1-x)^{\alpha}\frac{\ell(ry)}{\ell(ry(1-x))}\le
(1+\varepsilon)(1-x)^{\alpha-\varepsilon}\,.
\end{equation*}
We have
\begin{multline}
  \label{eq:79}
  \sup_{y\in(A/r,1]}\int_{1-\eta}^1
\ind_{\{ry(1-x)>A'\}}\,xy\,\frac{1- F (ry)}{(1-x)(1- F (ry(1-x)))}
\frac{ F (ry\,dx)}{1- F (r)}\\
\le(1+\varepsilon)\,\sup_{y\in(A/r,1]}
\frac{y(1- F (ry))}{1- F (r)}\,\sup_{z>A}
\frac{1}{1- F (z)}   \int_{1-\eta}^{1-A'/z} (1-x)^{\alpha-\varepsilon-1}
\, F (z\,dx)\,.
\end{multline}
By the first assertion of part 2 of Lemma~\ref{le:slowly},
 the first supremum on the right of \eqref{eq:79} is bounded from above
uniformly in all $A$ (and $r$) large enough.

For the second supremum, the integration by parts yields
\begin{multline*}
\frac{1}{1- F (z)}  \int_{1-\eta}^{1-A'/z}(1-x)^{\alpha-\varepsilon-1}
 F (z\,dx)=
(1+\varepsilon-\alpha)\int_{1-\eta}^{1-A'/z}(1-x)^{\alpha-\varepsilon-2}
\,\frac{ F (z)- F (zx)}{1- F (z)}\,dx\\
+\eta^{\alpha-\varepsilon-1}\,\frac{ F (z-A')- F (z(1-\eta))}{1- F (z)}\,.
\end{multline*}
By condition 2 of Theorem~\ref{the:conv_levy}, we have that, for $z$ large
enough, for $\eta$ small enough,
and for $x\in[1-\eta,1-A'/z]$,
\begin{equation*}
\frac{ F (z)- F (zx)}{1- F (z)}\le C(1-x)\,.
\end{equation*}
Therefore,
\begin{equation*}
\int_{1-\eta}^{1-A'/z}(1-x)^{\alpha-2-\varepsilon}
\,\frac{ F (z)- F (zx)}{1- F (z)}\,dx\le C
\int_{1-\eta}^1(1-x)^{\alpha-\varepsilon-1}\,dx\,.
\end{equation*}
It follows that
\begin{equation*}
    \lim_{\eta\to0}\lim_{A\to\infty}\sup_{z>A}\frac{1}{1- F (z)}
\int_{1-\eta}^1 (1-x)^{\alpha-\varepsilon-1}
\, F (z\,dx)=0\,.
\end{equation*}
By \eqref{eq:79}, we conclude that the third term on the rightmost side of
\eqref{eq:39} tends to zero as $A\to\infty$ and $\eta\to0$\,.
By letting successively $r\to\infty$, $A\to\infty$, and
$\eta\to0$ in \eqref{eq:39} and picking large enough $A'$, we obtain \eqref{eq:31}\,.

Now consider $I_3$\,. We have, by a change of variables,
\begin{eqnarray*}
  I_3 &\le &
M\int_{h(t,r)}^t \int_0^{1-\eta}
\bl(1-\widehat{X}^{(r)}(s)\br)x\,\ind_{\{\widehat{X}^{(r)}(s)<1\}}\,
(1-x)^{\alpha-1}
\,\frac{ F \bl(r(1-\widehat{X}^{(r)}(s))\,dx\br)}{1- F (r)}\,d s\\
\\
&\le &
M(t-h(t,r))\,
\frac{\eta^{\alpha-1}}{1- F (r)}\,
\int_0^{1-\eta}x F (r\,dx)\,.
\end{eqnarray*}
By the second assertion of part 2 of Lemma~\ref{le:slowly},
\begin{equation}
  \label{eq:81}
\limsup_{r\to\infty}  \frac{1}{1- F (r)}\,
\int_0^{1-\eta}x F (r\,dx)\le \frac{1}{1-\alpha}\,.
\end{equation}
Since
$t-h(t,r)\to0$ as
$r\to\infty$, it follows that
\begin{equation}
  \label{eq:36}
\lim_{r\to\infty}\sup_{\omega\in\Omega}I_3=0\,.
\end{equation}

For $I_4$, we reuse the earlier arguments to obtain that, for $A>0$,
\begin{multline}
  \label{eq:41}
    I_4\le
Mt\frac{A}{r(1- F (r))}
 \Big(\,\frac{1}{\displaystyle\inf_{z\le A}\mathbf{E}^\ast \exp(-\gamma\chi^{(z)})}+
\eta^{\alpha-1}\Big) F (A)\\
+
Mt\eta^{-1}
\,\sup_{z>A}\sup_{x\in[0,1-\eta]}\Big|\frac{(1-x)\mathbf{E}^\ast
\exp(-\gamma\chi^{(z(1-x))})}
{\mathbf{E}^\ast \exp(-\gamma\chi^{(z)})}-(1-x)^{\alpha}\Big|
\sup_{y\in[0,1]}
  \frac{y}{1- F (r)}\,
\int_0^{1-\eta}x\,
\, F (ry\,dx)\,.
\end{multline}
Take an arbitrary $\varepsilon\in(0,C_0)$\,.
By Lemma \ref{le:denomin}, $\mathbf{E}^\ast
\exp(-\gamma\chi^{(z)})\sim C_0/(z(1- F (z)))$ as $z\to\infty$. Therefore, 
for all $z$ large enough and uniformly over $x\in[0,1-\eta]$, we have
\begin{equation*}
  \Big|\frac{(1-x)\mathbf{E}^\ast\exp(-\gamma\chi^{(z(1-x))})}
{\mathbf{E}^\ast\exp(-\gamma\chi^{(z)})}-(1-x)^{\alpha}\Big|\le
\varepsilon\,\frac{1- F (z)}
{1- F (z(1-x))}+
\Big|\frac{1- F (z)}
{1- F (z(1-x))}-(1-x)^{\alpha}\Big|\,.
\end{equation*}
The first summand on the right  is not greater than
$\varepsilon$\,. The second summand tends to zero as $z\to\infty$
uniformly over $x\in[0,1-\eta]$ by the Uniform Convergence Theorem for
slowly varying functions.  
Thus,
\begin{equation*}
  \lim_{z\to\infty}\sup_{x\in[0,1-\eta]}
\Big|\frac{(1-x)\mathbf{E}^\ast\exp(-\gamma\chi^{(z(1-x))})}
{\mathbf{E}^\ast\exp(-\gamma\chi^{(z)})}-(1-x)^{\alpha}\Big|=0\,.
\end{equation*}
Since also \eqref{eq:81} holds, the second summand
on the right of \eqref{eq:41} tends to zero as $A$ and $r$ tend to infinity.
Since the first summand  tends to zero as $r\to\infty$, we arrive at the convergence
\begin{equation}
  \label{eq:37}
\lim_{r\to\infty}\sup_{\omega\in\Omega}I_4=0\,.
\end{equation}

We now consider $I_5$\,.
For $\delta\in(0,1-\eta)$, we employ a change of variables in
the inside integral of the first term on the left-hand side of the first inequality
and get
\begin{multline}
  \label{eq:44}
  I_5 \le
\int_0^t \int_0^{1-\eta}
\abs{V_3(\widehat{X}^{(r)}(s),  x)}
\,\frac{ F \bl(r(1-\widehat{X}^{(r)}(s))\,dx\br)}{1- F (r)}
\ind_{\{\widehat{X}^{(r)}(s)\in(1-\delta,1)\}}\,
d s\\
+
\int_0^t \int_0^{1-\eta}
\abs{V_2(\widehat{X}^{(r)}(s), x)}
\ind_{\{\widehat{X}^{(r)}(s)\in(1-\delta,1)\}}\,dx\,d s\\
+
\Big|\int_0^t \int_0^{1-\eta}
V_3(\widehat{X}^{(r)}(s), x)
\,\frac{ F \bl(r(1-\widehat{X}^{(r)}(s))\,dx\br)}{1- F (r)} \ind_{\{\widehat{X}^{(r)}(s)\le 1-\delta\}}
\,d s
\\-
 \int_0^t \int_0^{1-\eta}
V_2(\widehat{X}^{(r)}(s),x)\,\ind_{\{\widehat{X}^{(r)}(s)\le 1-\delta\}}
\,dx\,d s\Big|
\\
\le  Mt\eta^{\alpha-1}\,
\Bigg(\int_0^{\delta(1-\eta)}x\,\frac{ F (r\,dx)}{1- F (r)}
+\delta^{1-\alpha}\int_0^{1-\eta}\alpha
x^{-\alpha}\,dx\Bigg)
\quad \quad \quad \quad \quad \quad \quad \quad \quad \quad \quad \quad
\quad \quad
\\
+ t\sup_{y\in[\delta,1]}\Big| \int_0^{1-\eta}
g(yx)(1-x)^{\alpha-1}
\,\frac{ F \bl(ry\,dx\br)}{1- F (r)}-
 \int_0^{1-\eta}
g(yx)y^{-\alpha}(1-x)^{\alpha-1}
\alpha x^{-\alpha-1}\,dx\Big|\,. \quad \quad
\end{multline}
By  the second assertion of part 2 of Lemma~\ref{le:slowly},
we conclude that
 the term in the parentheses  on the rightmost side
of \eqref{eq:44} tends to zero as
$r\to\infty$ and $\delta\to0$\,, i.e.,
\begin{equation}
  \label{eq:88}
\lim_{\delta\to0}\limsup_{r\to\infty}
Mt\eta^{\alpha-1}\,
\Bigg(\int_0^{\delta(1-\eta)}x\,\frac{ F (r\,dx)}{1- F (r)}
+\delta^{1-\alpha}\int_0^{1-\eta}\alpha
x^{-\alpha}\,dx\Bigg)=0\,.
\end{equation}
 The other summand is not greater
than
\begin{multline}
  \label{eq:42}
 t\sup_{y\in[\delta,1]} \int_0^{1-\eta}
|g(yx)|\Big|\frac{1- F (ry)}{1- F (r)}-y^{-\alpha}\Big|(1-x)^{\alpha-1}
\,\frac{ F \bl(ry\,dx\br)}{1- F (ry)}+
   t\sup_{z\ge r\delta}\sup_{y\in[\delta,1]}
   W(0,y,z)
\\\le
M t\eta^{\alpha-1}\sup_{y\in[\delta,1]}\Big|\frac{1- F (ry)}{1- F (r)}
-y^{-\alpha}\Big|\sup_{y\in[\delta,1]} \int_0^{1-\eta}
yx\,\frac{ F \bl(ry\,dx\br)}{1- F (ry)}\\
+Mt\eta^{\alpha-1}\sup_{z\ge r\delta}\int_0^\delta x
\,\frac{ F \bl(z\,dx\br)}{1- F (z)}
+Mt\eta^{\alpha-1}\int_0^\delta \alpha x^{-\alpha}\,dx\,
+   t\sup_{z\ge r\delta}\sup_{y\in[\delta,1]}
W(\delta,y,z)
\end{multline}
where, for $0\le a \le 1$,
$$
W(a,y,z)=
\Bigg| \int_a^{1-\eta}
g(yx)y^{-\alpha}(1-x)^{\alpha-1}
\,\frac{ F \bl(z\,dx\br)}{1- F (z)}-
 \int_a^{1-\eta}
g(yx)y^{-\alpha}(1-x)^{\alpha-1}
\alpha x^{-\alpha-1}\,dx\Bigg| \,.
$$
Since $1- F (y)$ is regularly varying with index $-\alpha$,
\begin{equation*}
\sup_{y\in[\delta,1]}\Big|\frac{1- F (ry)}{1- F (r)}-y^{-\alpha}\Big|\to
0 \text{ as }r\to\infty.
\end{equation*}
Also,
\begin{equation*}
  \sup_{y\in[\delta,1]} \int_0^{1-\eta}
yx\,\frac{ F \bl(ry\,dx\br)}{1- F (ry)}\le
 \int_0^{1-\eta}
x\,\frac{ F \bl(r\,dx\br)}{1- F (r)}\,.
\end{equation*}
By the second assertion of
part 2 of Lemma~\ref{le:slowly}, we get 
\begin{equation}
  \label{eq:84}
  \lim_{r\to\infty}M t\eta^{\alpha-1}\sup_{y\in[\delta,1]}
\Big|\frac{1- F (ry)}{1- F (r)}
-y^{-\alpha}\Big|\sup_{y\in[\delta,1]} \int_0^{1-\eta}
yx\,\frac{ F \bl(ry\,dx\br)}{1- F (ry)}=0\,.
\end{equation}
Another application of
 the second assertion of part 2 of Lemma~\ref{le:slowly} yields
\begin{equation*}
\limsup_{z\to\infty}  \int_0^\delta x
\,\frac{ F \bl(z\,dx\br)}{1- F (z)}
\le \frac{\delta^{1-\alpha}}{1-\alpha}\,,
\end{equation*}
so
\begin{equation}
  \label{eq:85}
\lim_{\delta\to0}\lim_{r\to\infty}\sup_{z\ge r\delta}
  \int_0^\delta x
\,\frac{ F \bl(z\,dx\br)}{1- F (z)}=0\,.
\end{equation}
Also, clearly, 
\begin{equation}
  \label{eq:86}
  \lim_{\delta\to0}Mt\eta^{\alpha-1}\int_0^\delta \alpha x^{-\alpha}\,dx=0\,.
\end{equation}

Now consider the last term on the right-hand side of
\eqref{eq:42}\,.
Since $1- F $ is regularly varying at infinity with index $-\alpha$,
the  $ F \bl(z\,dx\br)/(1- F (z))$,
considered as measures on $[\delta,1-\eta]$,
 weakly converge as
$z\to\infty$ to  $ \alpha x^{-\alpha-1}\,dx$\,. Since
$g(x)$ is a continuous function, the  functions
$(g(yx)y^{-\alpha}(1-x)^{\alpha-1},\,x\in[\delta,1-\eta])$ are uniformly bounded
and equicontinuous
over $y\in[\delta,1]$\,. Therefore,
\begin{equation}
  \label{eq:87}
\lim_{z\to\infty}    \sup_{y\in[\delta,1]}
    W(\delta,y.z)
=0\,.
\end{equation}
Thus, by \eqref{eq:44}, \eqref{eq:88}, \eqref{eq:42}, \eqref{eq:84},
\eqref{eq:85}, \eqref{eq:86}, and \eqref{eq:87},
\begin{equation}
  \label{eq:38}
\lim_{r\to\infty}\sup_{\omega\in\Omega}I_5=0\,.
\end{equation}
Putting together \eqref{eq:30}, \eqref{eq:32}, \eqref{eq:31},
\eqref{eq:36}, \eqref{eq:37}, and \eqref{eq:38}, we conclude that
\begin{equation*}
\lim_{r\to\infty}
\sup_{\omega\in\Omega}\Big|\int_0^t\int_\R g(x)\,\widetilde{\nu}_3^{(r)}(ds,dx)-
\int_0^t\int_\R g(x)\,K(\widehat{X}^{(r)}(s);dx)\,ds\Big|=0\,.
\end{equation*}
On recalling \eqref{eq:29} and \eqref{eq:33}, we arrive at \eqref{eq:28}.

The assertion of the lemma for   nonlattice $F$  has
been proved. If $F$ is  lattice  with span $h$, the proof
proceeds analogously provided one assumes  that $r=nh$,
where $n\in\N$, that the y's are of
the form $k/n$ for $k=1,2,\ldots,n$, that the z's are of the form $kh$
for $k\in\N$, and that $n\to\infty$\,.
\end{proof}
\begin{proof}[Proof of Theorem~\ref{the:conv_levy}]
We apply Lemma~\ref{le:conv_jump} to the processes $\widehat{X}^{(r)}$
under the measures $\mathbf{Q}^{(r)}$ if $F$ is a nonlattice
distribution and to the processes $\widehat{X}^{(nh)}$
under the measures $\mathbf{Q}^{(nh)}$ if $F$ is a lattice
distribution with span $h$\,.
Conditions 1 and 2 of the lemma follow by
Lemma~\ref{le:function_K}. Condition 3 holds by
Lemma~\ref{le:conv_char}. Condition 4 is obviously met. The validity
of condition 5
had been established earlier.
\end{proof}

{\bf Acknowledgements.\/} This research was initiated during the second
author's visit to the Heriot-Watt University. The warm welcome of the Department of the
 Actuarial  Mathematics and Statistics at the Heriot-Watt University  and
the support of the European Commission
under the Marie Curie International Incoming Fellowship Programme
 are gratefully acknowledged.

The authors would like to thank the anonymous referees for their helpful comments
and suggestions. 

\appendix

\section{Auxiliary results}

\subsection{A proof of two results by Korshunov \cite{Kor05}}

The result formulated in this Subsection is contained in  Theorems 1 and 2
in Korshunov \cite{Kor05}. However, some details of the proof are
omitted there (especially for the lattice case), so we fill in the gaps
in our proof below.
We recall that $\mathbf{E}^\ast$ denotes expectation with respect to measure
$\mathbf{P}^\ast$ defined by \eqref{eq:48}.
  \begin{lemma}
    \label{le:denomin}
Let condition 1  of Theorem~\ref{the:conv_levy} hold.
\begin{enumerate}
\item If, in addition, $F$ is   nonlattice, then, for some
$C_0>0$, as $r\to\infty$\,,
    \begin{equation}\label{Dima2}
      \bE^\ast e^{-\gamma\chi^{(r)}}
\sim  \frac{C_0}{r(1- F (r))}\,.
    \end{equation}
\item If, instead, $F$ is  lattice with span $h$, then, for
  some $C_0'>0$, as
  $n\to\infty$,
    \begin{equation*}
      \bE^\ast e^{-\gamma\chi^{(nh)}}
\sim  \frac{C_0'}{nh(1- F (nh))}\,.
    \end{equation*}
\end{enumerate}

   \end{lemma}
   \begin{remark}
Note that, by Karamata's theorem, the coefficients $C_3$ in
\eqref{Dima} and $C_0$ in \eqref{Dima2} are related by
$C_0=C_3(1-\alpha)/\gamma$.
   \end{remark}
  \begin{proof}
    We introduce strict ascending  ladder indices $T_1$, $T_2$, ... by letting
$T_0=0$ and $T_n=\min\{k> T_{n-1}:\,S_k-S_{T_{n-1}}> 0\}$ for $n\in\N$.
Let $\zeta_k=S_{T_k}-S_{T_{k-1}}$ for $k\in\N$. Under
$\mathbf{P}^\ast$,
the $\zeta_k$ are a.s. finite
    and i.i.d., and $\mathbf{E}^\ast T_1<\infty$, see Asmussen
    \cite[VIII.2]{Asm03}.
 We let   $F_+$ denote  the common distribution function of the $\zeta_k$
    (under $\mathbf{P}^\ast$)\,.
Adapting the argument of the proof of  Lemma 2 in Korshunov
     \cite{Kor05}, we write, for $x\ge 0$,
     \begin{equation*}
       \frac{1-F_+(x)}{1-F(x)}=\int_{-\infty}^0 \frac{1-F(x-u)}{1-F(x)}\,H(du),
     \end{equation*}
where $H(u)=\ind_{\{u=0\}}+\sum_{k=1}^\infty \mathbf{P}^\ast(S_1\le 0,S_2\le
0,\ldots,S_{k}\le 0,S_{k}\le u)$ for $u\le 0$\,. Under condition 1 of
Theorem~\ref{the:conv_levy}, $\lim_{x\to\infty}(1-F(x-u))/(1-F(x))=1$, so by Lebesgue's
bounded convergence theorem,
\begin{equation*}
  \lim_{x\to\infty}\frac{1-F_+(x)}{1-F(x)}=H(0)\,.
\end{equation*}
Since $H(0)=1+\sum_{k=1}^\infty \mathbf{P}^\ast(T_1>k)
=\mathbf{E}^\ast T_1$, we conclude that
     \begin{equation}
       \label{eq:23}
       1-F_+(x)\sim  (1-F(x))\mathbf{E}^\ast T_1\;\;\text{ as }\;x\to\infty\,.
     \end{equation}
Thus, $1-F_+(x)$ is  regularly varying at infinity with index
$-\alpha$\,.

Since  $\chi^{(r)}$ is the overshoot over level $r$ of the random walk
$S_n$, it is also the overshoot over  $r$ of the random walk
 with steps $\zeta_k$.
Denoting by $H_+(x)$  the corresponding renewal function, i.e.,
$H_+(x)=\ind_{\{x\ge0\}}+\sum_{n=1}^\infty \bP^\ast(\sum_{k=1}^n
 \zeta_k\le x)$,
we have
\begin{equation*}
          \bE^\ast e^{-\gamma\chi^{(r)}}=
\int_{[0,r)} \int_{[r-x,\infty)} e^{-\gamma(y-(r-x))}\,F_+(dy)\,H_+(dx).
\end{equation*}
On introducing
\begin{equation}
      \label{eq:40}
  z(x)=\int_{[x,\infty)} e^{-\gamma(y-x)}\,F_+(dy)\,,
\end{equation}
we obtain that
\begin{equation}
  \label{eq:27}
    \bE^\ast e^{-\gamma\chi^{(r)}}=
\int_{[0,r]} z(r-x)\,H_+(dx)-z(0)\Delta H_+(r).
\end{equation}
Note that $z(x)=O(1/x)$ as $x\to\infty$, which follows from the
following calculations:
\begin{align*}
  z(x)=\int_{[x,x+\ln x/\gamma]}+\int_{(x+\ln
  x/\gamma,\infty)}&\le
(F_+(x+\frac{\ln x}{\gamma})-F_+(x-1))+\frac{1}{x}\,,\\
F_+(x+\frac{\ln x}{\gamma})-F_+(x-1)&
\sim \frac{\alpha}{\gamma}\, \frac{(1-F_+(x))\ln x}{x}\,\,,
\end{align*}
where for the  equivalence we used the fact that
$(1-F_+(x+\ln x/\gamma))/(1-F_+(x-1))\sim
(1+\ln x/(\gamma x))^{-\alpha}$ by the uniform convergence theorem for
regularly varying functions, see Bingham, Goldie, and Teugels
\cite[Theorem 1.5.2]{Reg89}\,.

Suppose now that $F$ is  nonlattice. Then $F_+$ is 
nonlattice  too, see Asmussen \cite[VIII.1]{Asm03}.
By \eqref{eq:40}, the function $z(x)$ is directly Riemann integrable, as defined in Feller
\cite[XI.1]{Fel71},  and
\begin{equation*}
    \int_0^\infty z(x)\,dx=\int_0^\infty e^{-\gamma x}(1-F_+(x))\,dx.
\end{equation*}
Thus, by Theorem 3 of Erickson \cite{Eri70},
 as $r\to\infty$,
\begin{equation*}
\int_{[0,r]} z(r-x)\,H_+(dx)\sim
\Bl(\int_0^r(1-F_+(x))\,dx\Br)^{-1}
\frac{\sin \pi\alpha}{\pi(1-\alpha)}\,
\int_0^\infty e^{-\gamma x}(1-F_+(x))\,dx\,.
\end{equation*}
In addition, by Theorem 1 of Erickson \cite{Eri70}, as $r\to\infty$,
\begin{equation*}
  \Delta H_+(r)\int_0^r(1-F_+(x))\,dx\to0\,.
\end{equation*}
If we also recall \eqref{eq:23} and the fact that, according to Karamata's
theorem (see Proposition 1.5.8 in Bingham, Goldie, and Teugels \cite{Reg89}),
$\int_0^r(1-F_+(x))\,dx\sim r(1-F_+(r))/(1-\alpha)$, we obtain the asymptotic
 equivalence asserted in part 1 with
\begin{equation*}
  C_0=\frac{1}{\mathbf{E}^\ast T_1}\,\frac{\sin\pi\alpha}{\pi}
\int_0^\infty e^{-\gamma x}(1-F_+(x))\,dx\,.
\end{equation*}
For lattice distributions, we haven't been able to find in the literature an analogue
of Erickson's Theorem 3. Therefore, we, in effect, deduce it from the
local renewal theorem of Garsia and Lamperti \cite{GarLam62} for our
particular case by using the approach of Erickson \cite{Eri70}.
As a matter of fact, we improve on Erickson's argument so that we can
give a streamlined proof of his Theorem 3.

Let $F$ be  lattice  with span $h$. Then
$F_+$ is also  lattice  with span $h$.
We can write for $\theta\in(0,1)$ and  suitable $A>0$, on recalling
that $z(x)=O(1/x)$ as $x\to\infty$,
\begin{multline}
  \label{eq:71}
  \int_{[0,nh]} z(nh-x)\,H_+(dx)=\int_{[0,\theta nh]} z(nh-x)\,H_+(dx)+
\int_{(\theta nh,nh]} z(nh-x)\,H_+(dx)\\\le
\frac{A}{(1-\theta)nh}\,H_+(\theta nh)+
\int_{(\theta nh,nh]} z(nh-x)\,H_+(dx)\,.
\end{multline}
By the fact that the tail of $F_+$ is regularly varying at infinity with index
$-\alpha$, we have (see Feller \cite[XIV.3]{Fel71} or Bingham, Goldie,
and Teugels \cite[8.6]{Reg89}) that
$  H_+(x)\sim (\sin \pi\alpha/\pi\alpha)(1-F_+(x))^{-1}$
as $x\to\infty$, so
\begin{equation}
  \label{eq:83}
\lim_{\theta\to0}\limsup_{n\to\infty} nh(1-F_+(nh))
  \frac{A}{(1-\theta)nh}\,H_+(\theta nh)=0\,.
\end{equation}
Next,
\begin{multline}
  \label{eq:89}
 nh(1-F_+(nh)) \int_{(\theta nh,nh]} z(nh-x)\,H_+(dx)=
nh(1-F_+(nh))\sum_{k=\lfloor \theta n\rfloor
}^{n}z\bl((n-k)h\br)\Delta H_+(kh)\\=
h\sum_{k=0 }^{\infty}z\bl(kh\br)\,\frac{n}{n-k}\,\frac{1-F_+(nh)}{1-F_+((n-k)h)}
\,m(n-k)\,\ind_{\{n-k\ge \lfloor \theta n\rfloor\}}\,,
\end{multline}
where we used the notation $m(k)=k(1-F_+(kh))\Delta H_+(kh)$ for
$k>0$ and $m(k)=0$ for $k\le0$.  Since the $\zeta_i$ assume
 values $kh,\,k\in\N$, it follows,
by Garsia and Lamperti \cite{GarLam62}, that
 \begin{equation}
   \label{eq:72}
\lim_{k\to\infty}m(k)=\frac{\sin\pi\alpha}{\pi}\,.
 \end{equation}
 Hence,
\begin{equation}
  \label{eq:90}
    \lim_{n\to\infty}\frac{n}{n-k}\,\frac{1-F_+(nh)}{1-F_+((n-k)h)}
\,m(n-k)\,\ind_{\{n-k\ge \lfloor \theta n\rfloor\}}=\frac{\sin\pi\alpha}{\pi}\,.
\end{equation}
Also by the uniform convergence theorem for regularly varying
functions,
\begin{equation*}
  \lim_{n\to\infty}\sup_{k\le n- \lfloor \theta n\rfloor}
\Big|\frac{1-F_+(nh)}{1-F_+((n-k)h)}-\frac{n^{-\alpha}}{(n-k)^{-\alpha}}\Big|=0\,.
\end{equation*}
Thus,
\begin{equation*}
  \limsup_{n\to\infty}\sup_{k=0,1,2,\ldots}
\frac{n}{n-k}\,\frac{1-F_+(nh)}{1-F_+((n-k)h)}
\,m(n-k)\,\ind_{\{n-k\ge \lfloor \theta n\rfloor\}}<\infty,
\end{equation*}
so by \eqref{eq:90}, Fatou's lemma, and Lebesgue's bounded convergence theorem,
\begin{equation}
  \label{eq:92}
    \lim_{n\to\infty}
\sum_{k=0 }^{\infty}z\bl(kh\br)\,\frac{n}{n-k}\,\frac{1-F_+(nh)}{1-F_+((n-k)h)}
\,m(n-k)\,\ind_{\{n-k\ge \lfloor \theta n\rfloor\}}=
\frac{\sin\pi\alpha}{\pi}\sum_{k=0 }^{\infty}z\bl(kh\br)\,.
\end{equation}
Putting together \eqref{eq:23}, \eqref{eq:71}, \eqref{eq:83},
\eqref{eq:89}, and \eqref{eq:92},
we conclude that
\begin{equation}
  \label{eq:91}
\lim_{n\to\infty}nh(1-F(nh))  \int_{[0,nh]} z(nh-x)\,H_+(dx)=
\frac{1}{\mathbf{E}^\ast T_1}\,\frac{\sin\pi\alpha}{\pi}\sum_{k=0 }^{\infty}z(kh)h\,.
\end{equation}
By Garsia and Lamperti \cite{GarLam62},
\eqref{eq:23}, and \eqref{eq:72}, $nh(1-F(nh))z(0)\Delta H_+(nh)\to
h\,z(0)(\mathbf{E}^\ast T_1)^{-1}\sin(\pi\alpha)/\pi$
as $n\to\infty$, so
the second assertion of the lemma follows by
\eqref{eq:27} and \eqref{eq:91} with
\begin{equation*}
  C_0'=\frac{1}{\mathbf{E}^\ast T_1}\,\frac{\sin\pi\alpha}{\pi}\,\sum_{k=1}^\infty z(kh)h=
\frac{1}{\mathbf{E}^\ast T_1}\,
\frac{\sin\pi\alpha}{\pi}\,\sum_{k=0}^\infty e^{-\gamma kh}\bl(1-F_+(kh)\br)h\,.
\end{equation*}
  \end{proof}

\subsection{Some properties of slowly and regularly varying functions}

The following  lemma comes in useful in the proof of
Theorem~\ref{the:conv_levy}.
Note that the first part is Potter's theorem (see Bingham, Goldie, and
Teugels \cite[Theorem  1.5.6]{Reg89}).
\begin{lemma}
  \label{le:slowly}
  \begin{enumerate}
  \item
Let $L(x)$ be a slowly varying at infinity  function. Then, given an
arbitrary $\varepsilon>0$, there exists $x_0>0$ such that
$L(x)/L(y)\le (1+\varepsilon)\bl((x/y)\vee(y/x)\br)^{\varepsilon}$ for all
$x\ge x_0$ and $y\ge x_0$\,.
\item
 If $ F $ is regularly varying at infinity with
index $-\alpha$, where $\alpha\in(0,1)$, then
\begin{equation*}
  \limsup_{\substack{r\to\infty\\A\to\infty}}
\sup_{y\in[A/r,1]}\frac{y(1- F (ry))}{1- F (r)}\le1\,
\end{equation*}
and,  for $y\in[0,1]$,
\begin{equation*}
\limsup_{r\to\infty}  \frac{1}{1- F (r)}\,
\int_0^{y}x F (r\,dx)\le \frac{y^{1-\alpha}}{1-\alpha}\,.
\end{equation*}
 \end{enumerate}
\end{lemma}
\begin{proof}
In order to prove the first inequality of part 2, note that the function
$\ell(x)=x^{\alpha}(1- F (x))$ is slowly varying at infinity. Hence, for
given arbitrary $\varepsilon\in(0,1-\alpha)$, we have by part 1 for all
$y\in(0,1]$ and $r$ such that $ry$ is large enough
\begin{equation*}
  \frac{y(1- F (ry))}{1- F (r)}=
\frac{y^{1-\alpha}\ell(ry)}{\ell(r)}\le
y^{1-\alpha}(1+\varepsilon)y^{-\varepsilon}\le 1+\varepsilon\,.
\end{equation*}
We prove the second inequality.
Integration by parts yields
\begin{equation*}
  \int_0^{y}
x\, F (r\,dx)=\int_0^y
\bl( F (ry)- F (rx)\br)\,dx\,.
\end{equation*}
On picking $A\in(0,ry)$ and
partitioning the integration interval $[0,y]$ into two
pieces $[0,A/r]$ and $(A/r,y]$,  we have
\begin{equation}
  \label{eq:20}
    \frac{1}{1- F (r)}\,
\int_0^y x F (r\,dx)\le\frac{A}{r(1- F (r))}+
  \frac{1}{1- F (r)}\,\int_{A/r}^y
\bl(1- F (rx)\br)\,dx\,.
\end{equation}
Let $\varepsilon\in(0,1-\alpha)$ be otherwise arbitrary.
If  $A$ is large enough, then by part 1, on
recalling that the function $x^\alpha(1- F (x))$ is slowly varying
at infinity and $y\le1$, we have
for all  $x\in[ A/r,y]$,
\begin{equation*}
  \frac{1- F (rx)}{1- F (r)}\le (1+\varepsilon)x^{-\alpha-\varepsilon}\,.
\end{equation*}
Therefore, for these  $A$ and $r$,
\begin{equation*}
    \frac{1}{1- F (r)}\,\int_{A/r}^y
\bl(1- F (rx)\br)\,dx \le
(1+\varepsilon)\int_0^y x^{-\alpha-\varepsilon}\,dx
= \frac{(1+\varepsilon)y^{1-\alpha-\varepsilon}}{1-\alpha-\varepsilon}\,.
\end{equation*}
The required bound follows now  by \eqref{eq:20} and the fact that
the first term on the right of \eqref{eq:20} tends to zero as
$r\to\infty$
(and as $A$ is kept fixed large enough)\,.
\end{proof}

\subsection{Convergence of the unconstrained random walk}
\label{sec:conv-cond-rand-1}

We recall that  $X$ is a stable subordinator   with L\'evy measure $ \alpha
 x^{-\alpha-1}\,dx\,,x>0,$ 
 and the processes $X^{(r)}$ are defined by
\eqref{eq:22}.
\begin{lemma}
  \label{le:conv_levy}
Let the right-hand tail of  $ F $ be regularly varying at infinity with
index $-\alpha$, where $\alpha\in(0,1)$\,.
Then the  $X^{(r)}$ under measure
$\mathbf{P}^\ast$  converge in distribution to  $X$\,.
\end{lemma}
The proof will use the following implication  of
 Lemma~\ref{le:conv_jump}. One can also apply
Theorem VII.3.4 on p.414 in Jacod and Shiryaev \cite{JacShi03} or Theorem 7.3.1
on p.592 in Liptser and Shiryaev \cite{lipshir}.
\begin{lemma}
  \label{le:conv-jump-trunc}
Consider a sequence $\Xi^{(n)}=(\Xi^{(n)}(t),\,t\in\R_+)$ of
$\R$-valued pure jump semimartingales of
locally bounded variation
with independent increments  defined on  filtered probability spaces
$(\Omega^{(n)},\mathcal{G}^{(n)},\mathbf{G}^{(n)},\mathbf{P}^{(n)})$. Let
$\Xi$ be a L\'evy 
process on a probability
space $(\Upsilon,\mathcal{G},\mathbf{\Pi})$, with L\'evy measure $K$ such that $\int_\R
1\wedge\abs{x}\,K(dx)<\infty$\,. Suppose that $\Xi^{(n)}(0)=0$
$\mathbf{P}^{(n)}$-a.s. and $\Xi(0)=0$ $\mathbf{\Pi}$-a.s. If for an arbitrary
$\R$-valued bounded continuous function $g(x),\,x\in\R,$
such that $\abs{g(x)}\le M\abs{x}$ in a neighbourhood of the
origin with some $M>0$,
\begin{equation*}
  \lim_{n\to\infty}\int_0^t\int_\R g(x)\,\nu^{(n)}(ds,dx)=
t\int_\R g(x)\,K(dx)\,,
\end{equation*}
where  $\nu^{(n)}(dt,dx)$ denotes the $\mathbf{G}^{(n)}$-predictable jump measure
 of $\Xi^{(n)}$\,, then the $\Xi^{(n)}$ converge in distribution to $\Xi$.
\end{lemma}
\begin{proof}[Proof of Lemma~\ref{le:conv_levy}]
The processes $X^{(r)}$ under $\mathbf{P}^\ast$ are pure jump
semimartingales of locally bounded variation with independent
increments. Their $\mathbf{F}^{(r)}$-predictable
 jump measures   under  $\bP^\ast$ are given by \eqref{eq:9}, so
for $g(x)$ as in the hypotheses,
\begin{equation}
  \label{eq:14}
  \int_0^t\int_\R g(x)\,\nu^{(r)}(ds,dx)=
 \lfloor (1- F (r))^{-1}\, t \rfloor
 \int_\R g(x) F (rdx)\,.
\end{equation}
We may assume that, for a suitable $M'$,
 $\abs{g(x)}\le M'\abs{x}$ for all $x$\,.
We have on writing $1- F (x)=x^{-\alpha}\ell(x)$, where $\ell$ is
 slowly varying at infinity,
\begin{equation*}
 \lfloor\, (1- F (r))^{-1}\, t \rfloor
 \Big|\int_{-\infty}^0 g(x) F (rdx)\Big|\le
M'\,\dfrac{r^{\alpha-1}}{\ell(r)}\, t \int_{-\infty}^0 \abs{x} F (dx)\,.
\end{equation*}
Since $\int_{-\infty}^0
\abs{x} F (dx)=\int_{-\infty}^0\abs{x}\exp(\gamma
x)\,\widehat{F}(dx)<\infty$ and $\alpha<1$, we obtain that
\begin{equation}
  \label{eq:68}
  \lim_{r\to\infty}
 \lfloor\, (1- F (r))^{-1}\, t \rfloor
 \int_{-\infty}^0 g(x) F (rdx)=0\,.
\end{equation}
For $\varepsilon>0$,
\begin{equation*}
 \lfloor\, (1- F (r))^{-1}\, t \rfloor
 \int_0^\varepsilon
\abs{  g(x)} F (rdx)\le
 t M'\int_0^\varepsilon x\,\frac{ F (rdx)}{1- F (r)}\,.
\end{equation*}
Hence, by the second assertion of part 2 of Lemma~\ref{le:slowly},
\begin{equation}
  \label{eq:19}
  \lim_{\varepsilon\to0}\limsup_{r\to\infty} \lfloor\, (1- F (r))^{-1}\, t \rfloor
 \int_0^\varepsilon   g(x) F (rdx)=0\,.
\end{equation}
The hypotheses on $ F $ imply that, for $x>\varepsilon>0$,
$( F (rx)- F (r\varepsilon))/(1- F (r))\to
\varepsilon^{-\alpha}-x^{-\alpha}$ as $r\to\infty$, so the $  F (dx)/(1- F (r))$, as
measures on $[\varepsilon,\infty)$, weakly converge to the measure $\alpha
x^{-\alpha-1}\,dx$\,. On recalling that $g(x)$ is a bounded and
continuous function, we conclude that
\begin{equation}
  \label{eq:57}
  \lim_{r\to\infty} \lfloor\, (1- F (r))^{-1}\,
  t \rfloor\int_\varepsilon^\infty
  g(x)\, F (dx)=t\int_\varepsilon^\infty g(x)\,\alpha x^{-\alpha-1}\,dx\,.
\end{equation}
By \eqref{eq:14}, \eqref{eq:68}, \eqref{eq:19}, and \eqref{eq:57},
\begin{equation*}
  \lim_{r\to\infty}\int_0^t\int_\R
  g(x)\,\nu^{(r)}(ds,dx)=t\int_0^\infty g(x)\,\alpha x^{-\alpha-1}\,dx\,,
\end{equation*}
which completes the proof by Lemma~\ref{le:conv-jump-trunc}.
\end{proof}

\def\cprime{$'$} \def\cprime{$'$} \def\cprime{$'$} \def\cprime{$'$}
  \def\polhk#1{\setbox0=\hbox{#1}{\ooalign{\hidewidth
  \lower1.5ex\hbox{`}\hidewidth\crcr\unhbox0}}} \def\cprime{$'$}
  \def\cprime{$'$} \def\cprime{$'$} \def\cprime{$'$}

\end{document}